\documentclass[12pt,twoside,a4paper]{amsart}
\usepackage{amssymb}
\usepackage{amsmath,amscd}
\usepackage[initials,nobysame,shortalphabetic]{amsrefs}

\def\dbar{\bar\partial}

\def\C{{\mathbb C}}

\def\Cn{\C^n}

\def\D{{\mathcal D}}
\def\M{{\mathcal M}}
\def\PM{{\mathcal{PM}}}

\def\Hom{{\rm Hom\, }}
\def\codim{{\rm codim\,}}

\def\E{{\mathcal E}}

\def\Ok{{\mathcal O}}
\def\Ow{{\tilde{\Ok}}}
\def\Oc{{\Ok_c}}

\def\Re{{\rm Re\,  }}

\def\reg{{\rm reg}}
\def\sing{{\rm sing}}
\newcommand{\Com}[1]{}
\DeclareMathOperator{\Id}{Id}
\DeclareMathOperator{\supp}{supp}
\DeclareMathOperator{\ann}{ann}

\def\be{\begin{equation}}
\def\ee{\end{equation}}

\newtheorem{thm}{Theorem}[section]
\newtheorem{lma}[thm]{Lemma}
\newtheorem{cor}[thm]{Corollary}
\newtheorem{prop}[thm]{Proposition}

\theoremstyle{definition}

\newtheorem{df}{Definition}

\theoremstyle{remark}

\newtheorem{preremark}{Remark}
\newtheorem{preex}{Example}

\newenvironment{remark}{\begin{preremark}}{\end{preremark}}
\newenvironment{ex}{\begin{preex}}{\end{preex}}

\numberwithin{equation}{section}

\begin{document}

\date{\today}

\address{R. L\"ark\"ang\\Department of Mathematics\\Chalmers University of Technology and the University of Gothenburg\\412 96 G\"OTEBORG\\SWEDEN}

\email{larkang@chalmers.se}




\author{Richard L\"ark\"ang}

\title[Residue currents and weakly holomorphic functions]{Residue currents associated with weakly holomorphic functions}

\begin{abstract}
    We construct Coleff-Herrera products and Bochner-Martinelli type residue currents associated with a tuple
    $f$ of weakly holomorphic functions, and show that these currents satisfy basic properties from the (strongly)
    holomorphic case. This include the transformation law, the Poincar\'e-Lelong formula and the equivalence of the Coleff-Herrera
    product and the Bochner-Martinelli type residue current associated with $f$ when $f$ defines a complete intersection.
\end{abstract}

\maketitle

\section{Introduction}

The basic example of a residue current, introduced by Coleff and Herrera in \cite{CH}, is a current
called the \emph{Coleff-Herrera product} associated with a strongly holomorphic mapping $f = (f_1,\dots,f_p)$.
The Coleff-Herrera product is defined by
\begin{equation} \label{chporig}
    \dbar \frac{1}{f_1}\wedge \dots \wedge \dbar \frac{1}{f_p} . \varphi =
    \lim_{\delta\to 0^+ }\int_{\cap \{ |f_j| = \epsilon_j(\delta) \} } \frac{\varphi}{f_1\dots f_p},
\end{equation}
where $\varphi$ is a test form and $\epsilon(\delta)$ tends to $0$ along a so-called \emph{admissible path},
which means essentially that $\epsilon_1(\delta)$ tends to $0$ much faster than $\epsilon_2(\delta)$ and so on,
for the precise definition, see \cite{CH}.
The Coleff-Herrera product was defined over an analytic space, however, most of the work on residue currents thereafter
has focused on the case of holomorphic functions on a complex manifold.
The theory of residue currents has various applications, for example to effective versions of division problems etc.,
see for example \cite{AW1}, \cite{BGVY}, \cite{TY} and the references therein.

On an analytic space $Z$,
with structure sheaf $\Ok_Z$,
the most common notion of holomorphic functions are the \emph{strongly holomorphic} functions, that is,
sections of the structure sheaf, or more concretely,
functions which are locally the restriction of holomorphic functions in any local embedding.
We will throughout the article assume that $Z$ is an analytic space of pure dimension.
In some cases, this can be a little too restrictive, and the \emph{weakly holomorphic} functions might be more
natural. These are functions defined on $Z_\reg$, which are holomorphic on $Z_\reg$ and locally bounded at $Z_\sing$.
Two reasons why these are natural:
the ring $\Ow_{Z,z}$ of germs of weakly holomorphic functions at $z$ is the integral
closure of $\Ok_{Z,z}$ in the ring $\M_{Z,z}$ of germs of meromorphic functions at $z$,
and weakly holomorphic functions correspond to strongly holomorphic functions in any normal modification of
$(Z,\Ok_Z)$.
A slightly better behaved but more restrictive notion are the \emph{c-holomorphic} functions,
denoted $\Ok_c$,
functions which are weakly holomorphic and continuous on all of $Z$.
We will throughout this article assume that $Z$ is an analytic space of pure dimension.

In a recent article \cite{Den}, Denkowski introduced a residue calculus for c-holomorphic functions,
and showed that this calculus satisfies many of the basic properties known from the strongly holomorphic
or smooth cases. It is then a natural question to ask what happens in the case of weakly holomorphic
functions. However, as in the c-holomorphic case, it is not obvious how to define the associated residue currents.

In the strongly holomorphic case, there are various ways to define the Coleff-Herrera product
(for the equivalence of various definitions of the Coleff-Herrera product, also in the non complete intersection case,
see for example \cite{LS}).
The definition we will use is based on analytic continuation as in \cite{Y}, which was inspired by the ideas in \cite{At}
and \cite{BG} that the principal value current $1/f$ of a holomorphic function $f$ can be defined by $(|f|^{2\lambda}/f) |_{\lambda = 0}$.
If $f = (f_1,\dots,f_p)$ is strongly holomorphic
on $Z$, we define the Coleff-Herrera product of $f$ by
\begin{equation*}
    \left.\frac{\dbar |f_1|^{2\lambda_1}\wedge \dots \wedge \dbar |f_p|^{2\lambda_p}}{f_1\dots f_p}\right\lvert_{\lambda_p = 0,\dots,\lambda_1 = 0},
\end{equation*}
where we by
$|_{\lambda_p = 0,\dots,\lambda_1 = 0}$
mean that we take the analytic continuation in
$\lambda_p$ to $\lambda_p = 0$,
then in
$\lambda_{p-1}$
and so on, see Section \ref{sectchp} for details.
Recall that a \emph{modification} of an analytic space $Z$ is a proper surjective holomorphic mapping $\pi : Y \to Z$ from an analytic space $Y$
such that there exists a nowhere dense analytic set $E \subset X$ with
$\pi|_{Y\setminus \pi^{-1}(E)} : Y \setminus \pi^{-1}(E) \to X\setminus E$ being a biholomorphism.
It is easy to see by analytic continuation, that if $\pi : Y \to Z$ is a modification of $Z$,
then the Coleff-Herrera product of $f$ can be defined as the push-forward
of the Coleff-Herrera product of $f' := \pi^{*}f$.
For weakly holomorphic functions, we can use this observation to define the Coleff-Herrera product, since the
pull-back of a weakly holomorphic function to the normalization is strongly holomorphic.
If $f$ is weakly holomorphic, we define the Coleff-Herrera product of $f$ by
\begin{equation} \label{eqchfirstdef}
    \mu^f := \dbar\frac{1}{f_1}\wedge \dots \wedge \dbar\frac{1}{f_p} := \pi_*\left( \dbar\frac{1}{f_1'}\wedge \dots \wedge \dbar\frac{1}{f_p'}\right),
\end{equation}
where $f' = \pi^* f$.
By the observation above, this of course coincides with the usual definition in case of strongly holomorphic functions,
and this definition is also consistent with the one in \cite{Den} in the case of c-holomorphic functions, see Proposition~\ref{denequal}.

Because of our definition, the properties we prove of the Coleff-Herrera product for weakly
holomorphic functions can mostly be reduced (by going back to the normalization) to the strongly holomorphic case.
Thus the main part of this article concerns giving a coherent exposition of the basic theory of residue currents in the
strongly holomorphic case. This is done based on analytic continuation of currents
and the notion of pseudomeromorphic currents as introduced in \cite{AW2}, which is developed
on a complex manifold. We will see that this approach works well also with strongly holomorphic
functions on an analytic space, and we believe that this might be of independent interest,
although most of the results should be known.

However, even for the statement of these properties in the weakly holomorphic case, two problems occur,
namely how is multiplication of a weakly holomorphic function with a current defined,
and what is the zero set of a tuple of weakly holomorphic functions? And hence also, what should a complete intersection mean?

With regards to defining multiplication of a weakly holomorphic function with a current,
we take a similar approach as for the definition of the Coleff-Herrera product.
Assume $\mu$ is a current on $Z$, and that there exists a modification $\pi : Y \to Z$,
with a current $\mu'$ on $Y$ such that $\mu = \pi_* \mu'$ (the existence of such $\mu'$ is
guaranteed if $\mu$ is pseudomeromorphic and $Y$ is the normalization of $Z$,
see the introduction of Section~\ref{sectmult}).
If $g$ is strongly holomorphic on $Z$, then
\begin{equation} \label{eqwholomult}
    g\mu = \pi_*( \pi^*g \mu').
\end{equation}
The right-hand side of \eqref{eqwholomult} still exists if $g$ is weakly holomorphic on $Z$
and $Y$ is normal, so we take this as a definition of $g\mu$.
However, that this is well-defined depend on the fact that we have a certain ``canonical'' representative
of the Coleff-Herrera product in the normalization (or any normal modification).
We will see in Section~\ref{sectmult} that \eqref{eqwholomult} depends on the choice of representative
$\mu'$ and can thus not be used to define a general multiplication of weakly holomorphic with
currents on $Z$.

For the zero set of one weakly holomorphic function, all reasonable definitions should coincide.
For the zero set of a weakly holomorphic mapping $f$, it is natural to take into account that
the zero sets of the individual components of $f$ can ``belong'' to different irreducible components.
We introduce in Section~\ref{sectzerosets} a notion of common zero set of $f$,
depending on $f$ as a mapping, and not only on the individual components,
which however may differ from the intersection of the respective zero sets.

The Coleff-Herrera product $\mu^f$ in \eqref{eqchfirstdef} associated with a
strongly holomorphic mapping $f = (f_1,\dots,f_p)$ satisfies
\begin{equation*}
    \supp \mu^f \subseteq Z_f \quad \text{and} \quad \dbar \mu^f = 0,
\end{equation*}
where $Z_f$ is the common zero set of $f$.
In addition, if $f$ forms a complete intersection, the Coleff-Herrera product is alternating
in the residue factors and
\begin{equation} \label{eqannmufinclusion}
    (f_1,\dots,f_p) \subseteq \ann \mu^f,
\end{equation}
where $(f_1,\dots,f_p)$ is the ideal generated by $f_1,\dots,f_p$, and $\ann \mu^f$ is the annihilator
of $\mu^f$, i.e., the ideal of holomorphic functions $g$ such that $g \mu^f = 0$.
We also have the \emph{transformation law} for residue currents (see \cite{DS}),
which says that if $f = (f_1,\dots,f_p)$ and $g = (g_1,\dots,g_p)$ define complete intersections,
and there exists a matrix $A$ of holomorphic functions such that $g = Af$, then
\begin{equation*}
    (\det A) \dbar \frac{1}{g_1}\wedge \dots \wedge \dbar \frac{1}{g_p} = \dbar \frac{1}{f_1}\wedge \dots \wedge \dbar \frac{1}{f_p}.
\end{equation*}
The \emph{Poincar\'e-Lelong formula} relates the Coleff-Herrera product of $f$
and the integration current $[Z_f]$ on $Z_f$ (with multiplicities) and it says that
\begin{equation*}
    \frac{1}{(2\pi i)^p}\dbar\frac{1}{f_1}\wedge \dots \wedge \dbar\frac{1}{f_p}\wedge df_p \wedge \dots \wedge df_1 = [Z_f].
\end{equation*}
We will see that in fact all those statements still hold also in the weakly holomorphic case.
However, as mentioned above, zero sets of weakly holomorphic functions and multiplication
of currents with weakly holomorphic functions need to be interpreted in the right way.

\begin{remark}
The inclusion \eqref{eqannmufinclusion} if $f$ defines a complete intersection
is one direction of the \emph{duality theorem} proven in \cite{DS} and \cite{P2},
which says that on a complex manifold, the inclusion is in fact (locally) an equality.
However, in \cite{L}, we show that on any singular variety, one can always find a tuple $f$
of strongly holomorphic functions such that the inclusion \eqref{eqannmufinclusion} is strict.
\end{remark}

Bochner-Martinelli type residue currents were first introduced in \cite{PTY} by Passare, Tsikh and Yger 
(on a complex manifold) as an alternative way of defining a residue current corresponding to a tuple of holomorphic functions.
In \cite{BVY}, Bochner-Martinelli type residue currents were constructed on an analytic space in order to
prove a generalization of Jacobi's residue formula, generalizing previous results in \cite{VY}
in the smooth case.

The Bochner-Martinelli type residue currents give another reason why our definition of Coleff-Herrera product is a natural one.
In the smooth case, it was proved in \cite{PTY} that if the functions define a complete intersection,
then the Coleff-Herrera product and the Bochner-Martinelli current coincide.
It is suggested in \cite{BVY} that the same statement holds in the singular case with a similar proof.
We will construct Bochner-Martinelli type residue currents associated with a tuple of weakly holomorphic functions,
and we will show that the equality between the Coleff-Herrera product and the Bochner-Martinelli type residue current holds
both in the strongly and weakly holomorphic cases.
An advantage of the Bochner-Martinelli
current, compared to the Coleff-Herrera product,
in the weakly holomorphic case is that it can be defined 
intrinsically on $Z$ as the analytic continuation of an arbitrarily smooth (depending on a parameter $\lambda$) form on $Z$.
In contrast, the Coleff-Herrera product is only defined as the analytic continuation of an arbitrarily smooth form on the
normalization of $Z$.

\section{Zero sets of weakly holomorphic functions} \label{sectzerosets}
    The behavior of the currents we define will depend in a crucial way on the zero sets of the weakly
    holomorphic functions, and in this section we will define the zero set of a weakly holomorphic mapping.
    
    \begin{df}
        Let $f \in \Ow(Z)$.
        If $f$ is not identically zero on all irreducible components of $Z$, we define the \emph{zero set}
        of $f$ by $Z_f := \{ z \in Z \ \lvert\ (1/f)_z \notin \Ow_z \}$.
        Let $Z_\alpha$ be the irreducible components of $Z$ where $f$
        is identically zero, and let $Z' = \overline{Z \setminus \cup_\alpha Z_\alpha}$.
        Then $f$ does not vanish identically on any of the irreducible components of $Z'$,
        and we define $Z_f$ as $\cup_\alpha Z_\alpha \cup Z_{f|_{Z'}}$.
    \end{df}

    \begin{remark} \label{zerocharacterization}
        We have $z \in Z_f$ if and only if there exists a sequence $z_i \to z$ with $z_i \in Z_\reg$ such that $f(z_i) \to 0$
        (since if we cannot find such a sequence, then $1/f$ is weakly holomorphic). Hence, when $f$ is c-holomorphic,
        $Z_f$ coincides with the usual zero set of $f$, when $f$ is seen as a continuous function.
    \end{remark}

    We will use the following characterization of the zero set of a weakly holomorphic function.
    However, since this is a special case of Proposition~\ref{zeroproperties}, we omit the proof.

    \begin{lma} \label{normalizationzeroset}
        Let $\pi : Z' \to Z$ be the normalization of $Z$.
        If $f \in \Ow(Z)$, then $Z_f$ is an analytic subset of $Z$, and $Z_f = \pi(Z_{\pi^* f})$. 
    \end{lma}

    We recall that an analytic space $Z$ is \emph{normal} if $\Ok_{Z,z} = \Ow_{Z,z}$ for all $z \in Z$,
    and that the \emph{normalization} $Z'$ of an analytic space $Z$ is the unique normal space $Z'$
    together with a proper finite surjective holomorphic mapping $\pi : Z' \to Z$ such that
    $\pi|_{Z' \setminus \pi^{-1}(Z_\sing)} : Z'\setminus \pi^{-1}(Z_\sing) \to Z_\reg$ is
    a biholomorphism, see for example \cite{G}.

    For any meromorphic function $\phi$, there is a standard notion of zero set of $\phi$,
    that we denote by $Z'_\phi$, which is defined by $Z'_\phi = \{ z \in Z\ \lvert\ (1/\phi)_z \notin \Ok_z \}$.
    Since weakly holomorphic functions are meromorphic, this gives another definition of zero set
    if $f$ is a weakly holomorphic function.
    Clearly $Z_f \subseteq Z'_f$, but as we see in the following example, the inclusion is in general strict,
    so the two definitions do not coincide.

    \begin{ex}
        Let $Z = \{ z^3 - w^2 = 0 \} \subseteq \C^2$, which has normalization $\pi(t) = (t^2,t^3)$,
        and let $f = 1 + w/z$. Since $\pi^* f = 1 + t^3/t^2 = 1 + t$, $f$ is weakly holomorphic on $Z$.
        Since $\{ \pi^* f = 0 \} = \{ t = -1 \}$, we get by Lemma~\ref{normalizationzeroset} that
        \begin{equation*}
        Z_f = \pi(\{ t = -1 \}) = \{ (1,-1) \}.
        \end{equation*}
        However,
        \begin{equation*}
            Z_f' = P_{1/f} = Z \cap \{ z+w = 0 \} = \{ (t^2,t^3)\ |\ t^2 = -t^3 \} = \{ (0,0), (1,-1) \},
        \end{equation*}
        so $Z_f \subsetneq Z_f'$.
    \end{ex}

    To study the dimension of zero sets of weakly holomorphic functions, 
    we will need the following lemma, which shows that subvarieties of the normalization
    correspond to subvarieties of $Z$ of the same dimension, and vice versa.

    \begin{lma} \label{dimnormalization}
        Let $\pi : Z' \to Z$ be the normalization of $Z$. If $Y'$ is a subvariety
        of $Z'$, then $\pi(Y')$ is a subvariety of $Z$ with $\dim Y' = \dim \pi(Y')$,
        and if $Y$ is a subvariety of $Z$,
        then $\pi^{-1}(Y)$ is a subvariety of $Z'$ with $\dim Y = \dim \pi^{-1}(Y)$.
    \end{lma}

    \begin{proof}
        The first part follows from Remmert's proper mapping theorem, when formulated as for example in \cite{G},
        since $\pi$ is a finite proper holomorphic mapping.
        We get from the first part that $\dim \pi^{-1}(Y) = \dim \pi(\pi^{-1}(Y)) = \dim Y$,
        where the second equality holds since $\pi$ is surjective.
    \end{proof}

    If $f \in \Ow(Z)$ and $f \not\equiv 0$ on any irreducible component of $Z$, then $\codim Z_f = 1$
    or $Z_f = \emptyset$. In fact, if $f' = \pi^* f$ and $Z_{f'} \neq \emptyset$, then $f'$ is strongly holomorphic, and $Z_{f'} = \{f' = 0\}$
    has codimension $1$, and since $Z_f = \pi(Z_{f'})$ by Lemma \ref{normalizationzeroset},
    $Z_f$ has codimension $1$ by Lemma~\ref{dimnormalization}.
    However, as is well-known, in contrast to the smooth case, subvarieties of codimension $1$ cannot in general be
    defined as the zero set of one single strongly holomorphic function. As we will see in the next example,
    this is the case in general for zero sets of weakly holomorphic functions, 
    even for c-holomorphic functions on an irreducible space.

    \begin{ex} \label{zerosetnotholo}
        Let $V = \{ z_1^3 - z_2^2 = z_3^3 - z_4^2 = 0 \} \subset \C^4$. Then $V$ has normalization
        $\pi : \C^2 \to V$, $\pi(t_1,t_2) = (t_1^2,t_1^3,t_2^2,t_2^3)$, and hence $f = z_2/z_1 - z_4/z_3$ is c-holomorphic
        since $\pi^* f = t_1 - t_2$. The set $Z_f = \{ (t^2,t^3,t^2,t^3) \}$ has codimension $1$ in $Z$.
        However, there does not exist a holomorphic function in a neighborhood of $0$ such that
        $f(t_1^2,t_1^3,t_2^2,t_2^3) = 0$ exactly when $t_1 = t_2$, since in that case, we could write
        $f(t_1^2,t_1^3,t_2^2,t_2^3) = (t_1-t_2)^m u(t_1,t_2)$ for some $m \in \mathbb{N}$, where $u(0,0) \neq 0$,
        which is easily seen to be impossible.
        Hence, $Z_f$ is not the zero set of one single strongly holomorphic function.
    \end{ex}

    \begin{ex} \label{intersectweakly}
        Let $Z = Z_1\cup Z_2 \subset \C^6$, where $Z_1 = \C^3\times \{0\}$ and $Z_2 = \{0\}\times\C^3$.
        Define the functions $f$ and $g$ by
        \begin{equation*}
            f(z) = \left\{ \begin{array}{cc} z_1 & z \in Z_1\setminus\{0\} \\ 1 & z \in Z_2\setminus\{0\} \end{array}\right. \text{and }
                g(z) = \left\{ \begin{array}{cc} 1 & z \in Z_1\setminus\{0\} \\ z_4 & z \in Z_2\setminus\{0\} \end{array}\right. .
        \end{equation*}
        Then $f,g \in \Ow(Z)$,
        and $Z_f = Z_1\cap \{ z_1 = 0\}$, and $Z_g = Z_2\cap \{ z_4 = 0 \}$ which both have codimension $1$
        in $Z$. However, $Z_f \cap Z_g = \{ 0 \}$, which has codimension $3$. Hence, zero sets of weakly holomorphic functions
        do not behave as well as one could hope with respect to intersections.
        If we let $f_1 = f_2 = f$, $f_3 = g$, then $Z_{f_1}\cap Z_{f_2} \cap Z_{f_3} = \{ 0 \}$ has codimension
        $3$, while $Z_{f_1} \cap Z_{f_2} = Z_f$ has codimension $1$ at $0$ in $Z$.
        Hence, if one defines a complete intersection for zero sets of weakly holomorphic functions
        $f = (f_1,\cdots,f_p)$ by requiring that $Z_{f_1}\cap\cdots\cap Z_{f_p}$ has codimension $p$ in $Z$,
        then it will not follow in general that $(Z_{f_1}\cap\cdots\cap Z_{f_k},z)$ has codimension $k$ for
        $z \in Z_{f_1}\cap\cdots\cap Z_{f_p}$.
    \end{ex}

    \begin{remark} \label{intersectcholo}
    Note that for c-holomorphic functions $f=(f_1,\cdots,f_p)$, if $f' = \pi^* f$, where $\pi : Z' \to Z$ is the normalization, then
    $\pi(Z_{f'_1}\cap\cdots\cap Z_{f'_p}) = Z_{f_1}\cap \cdots \cap Z_{f_p}$. Thus if we say that $f = (f_1,\cdots,f_p)$, where $f_i \in \Ok_c(Z)$,
    forms a complete intersection in $Z$ if $Z_{f_1}\cap\cdots\cap Z_{f_p}$ has codimension $p$, then this holds if and only if
    $f'$ forms a complete intersection in $Z'$ by Lemma~\ref{dimnormalization}.
    \end{remark}
    As we see in Example~\ref{intersectweakly}, this remark does not hold for weakly holomorphic functions,
    because there, $Z_{f}\cap Z_{g} = \{ 0 \}$,
    while $Z_{f'} \cap Z_{g'} = \emptyset$. Thus, the straightforward generalization of complete intersection, where the zero set
    $Z_{f_1}\cap\cdots\cap Z_{f_p}$ is required to have codimension $p$ does not share the same good properties in the weakly holomorphic case as
    in the strongly holomorphic (or c-holomorphic) case. Because of this, we will use a different definition of both the common zero set of
    weakly holomorphic functions and of a complete intersection. It coincides with the usual definitions in case of strongly holomorphic
    or c-holomorphic functions, and with our definition the problems above disappear.

    \begin{df} \label{defzerosetci}
        Let $f = (f_1,\cdots,f_p)$ be weakly holomorphic. We define the \emph{common zero set} of $f$, denoted by $Z_f$,
        as the set of $z \in Z$ such that there exists a sequence $z_i \in Z_\reg$ with $z_i \to z$, and $f_k(z_i) \to 0$
        for $k = 1,\cdots,p$. We will see that $Z_f$ is an analytic subset of $Z$, and hence we say that $f$ forms a
        \emph{complete intersection} if $Z_f$ has codimension $p$ in $Z$.
    \end{df}

    Note that by Remark~\ref{zerocharacterization}, this definition is consistent with the definition of $Z_f$ in the case of one function.
    We also see that in Example~\ref{intersectweakly}, $Z_{(f,g)} = \emptyset$, and hence, $(f,g)$ is not a complete intersection in our sense.
    Just as for one function, we can give a characterization of the zero set with the help of the normalization.
    \begin{prop} \label{zeroproperties}
        Let $f = (f_1,\cdots,f_p)$ be weakly holomorphic, and let $f' = \pi^* f$, where $\pi : Z' \to Z$ is the normalization.
        Then
        \begin{equation} \label{zfnormalization}
            Z_f = \pi(Z_{f'_1} \cap\cdots\cap Z_{f'_p}),
        \end{equation}
        and if $Z_f$ is nonempty, then it
        is an analytic subset of $Z$ of codimension $\leq p$. In general, 
        \begin{equation} \label{cholozeroset}
            Z_f \subseteq Z_{f_1}\cap\cdots\cap Z_{f_p},
        \end{equation}
        with equality if $f$ is c-holomorphic.
        In addition, $f$ is a complete intersection if and only if $f'$ is a complete
        intersection in the normalization.
    \end{prop}

    \begin{proof}
        If $z' \in Z_{f'_1}\cap\cdots\cap Z_{f'_p}$, then we can take a sequence $z_i' \to z'$ such that $z_i' \in \pi^{-1}(Z_\reg)$.
        Then, if we let $z_i = \pi(z_i')$, we get that $f_k(z_i) \to 0$, and hence we have the inclusion
        $Z_f \supseteq \pi(Z_{f'_1}\cap\cdots\cap Z_{f'_p})$ in \eqref{zfnormalization}.
        For the other inclusion, if we have a sequence $z_i \to z$ such that $z \in Z_f$, since $\pi$ is proper we can
        choose a convergent subsequence $z'_{k_i} \to z'$ such that $\pi(z'_{k_i}) = z_{k_i}$, and since $z \in Z_f$,
        we must have $f'(z') = 0$, so $z = \pi(z')$, with $z' \in Z_{f_1'}\cap\cdots\cap Z_{f_p'}$. 
        Now, the fact that $Z_f$ is an analytic subset of $Z$ follows by \eqref{zfnormalization} and Remmert's proper mapping
        theorem, since $Z_{f_i'}$ are analytic subsets of $Z'$.
        Since $f'$ is strongly
        holomorphic, $Z_{f'}$ has codimension $\leq p$, so by \eqref{zfnormalization} combined with Lemma~\ref{dimnormalization}
        we get that $Z_f$ has codimension $\leq p$.
        If $f$ is c-holomorphic, the equality in \eqref{cholozeroset} follows by \eqref{zfnormalization}
        since for any continuous mapping $f$, $Z_{f_1}\cap\cdots\cap Z_{f_p} = \pi(Z_{\pi^* f_1}\cap\cdots\cap Z_{\pi^* f_p})$,
        and the general case also follows from \eqref{zfnormalization} since $\pi(Z_{f_1'}\cap\cdots\cap Z_{f_p'}) \subseteq \pi(Z_{f_1'})\cap\cdots
        \cap \pi(Z_{f_p'}) = Z_{f_1}\cap\cdots\cap Z_{f_p}$.
        Finally, the fact that $f$ is a complete intersection if and only if $f'$ is a complete intersection follows from \eqref{zfnormalization}
        together with Lemma~\ref{dimnormalization}.
    \end{proof}

    We note that if $Z_{f_1}\cap\cdots\cap Z_{f_p}$ has codimension $\geq p$, then either $Z_f = \emptyset$, or $Z_f$ has codimension
    $p$ since by Proposition~\ref{zeroproperties}, $Z_f \subseteq Z_{f_1}\cap\cdots\cap Z_{f_p}$, and $Z_f$ has codimension at most $p$.
    Thus, we could have taken as definition of a complete intersection, that $Z_{f_1}\cap\cdots \cap Z_{f_p}$ has codimension $\geq p$,
    and our results about complete intersection would still be true.
    However, it would in general give weaker statements, since it since it might very well happen that
    $Z_{f_1}\cap\cdots\cap Z_{f_p}$ has codimension $< p$, while $Z_f$ has codimension $p$.
    In addition, results depending on the exact zero set, like the Poincar\'e-Lelong formula, Proposition~\ref{plthm},
    would of course not be true if one would use $Z_{f_1}\cap\dots\cap Z_{f_p}$ instead of $Z_f$.

    Note also that, if $f = (f_1,\cdots,f_p)$ is a complete intersection and $f_0 = (f_1,\cdots,f_k)$,
    then $(Z_{f_0},z)$  has codimension $k$ for $z \in Z_f$, since if $z' \in \pi^{-1}(z)$, then
    $(Z_{f_0'},z')$ has codimension $k$, and hence since $\pi$ is a finite proper holomorphic mapping,
    $(Z_{f_0},z) = \cup_{z_j'\in \pi^{-1}(z)} \pi( (Z_{f_0'},z_j') )$ has codimension $k$ in $Z$.

\section{Pseudomeromorphic currents on an analytic space} \label{sectpmcurrents}

    We will in this section introduce pseudomeromorphic currents on an analytic space. Pseudomeromorphic currents
    on a complex manifold were introduced by Andersson and Wulcan in \cite{AW2}, inspired by the fact that currents
    like the Coleff-Herrera product and Bochner-Martinelli type residue currents are pseudomeromorphic.
    Two important properties of pseudomeromorphic currents in the smooth case are the direct analogues of Proposition~\ref{pmexist}
    and Proposition~\ref{pseudom0}. Since these hold also in the singular case, many properties of residue currents
    hold also for strongly holomorphic functions by more or less the same argument as in the smooth case.

    The pseudomeromorphic currents are intrinsic objects of the analytic space $Z$, so we begin with explaining what
    we mean by a current on an analytic space. 
    We will follow the definitions used in \cite{BH} and \cite{HL}.
    To begin with, we assume that $Z$ is an analytic subvariety of $\Omega$, for some open set $\Omega \subseteq \Cn$.
    Then, we define the set of smooth forms of bidegree $(p,q)$ in $Z$ by $\E_{p,q}(Z) = \E_{p,q}(\Omega)/\mathcal{N}_{p,q,Z}(\Omega)$,
    where $\E_{p,q}(\Omega)$ are the smooth $(p,q)$-forms in $\Omega$ and $\mathcal{N}_{p,q,Z}(\Omega) \subset \E_{p,q}(\Omega)$ are
    the smooth forms $\varphi$ such that $i^*\varphi \equiv 0$, where $i : Z_\reg \to \Omega$ is the inclusion map.
    The set of test forms on $Z$, $\D_{p,q}(Z)$, are the forms in $\E_{p,q}(Z)$ with compact support.
    With the usual topology on $\D_{p,q}(\Omega)$ by uniform convergence of coefficients of differential forms together
    with their derivatives on compact sets, we give $\D_{p,q}(Z)$ the quotient topology from the projection
    $\D_{p,q}(\Omega) \to \D_{p,q}(Z)$. Then, $(p,q)$-currents on $Z$, denoted $\D'_{p,q}$, are the continuous linear functionals
    on $\D_{k-p,k-q}(Z)$, where $k = \dim Z$.
    However, more concretely, this just means that if $\mu$ is a $(p,q)$-current on $Z$, then $i_*\mu$ is a
    $(n-k + p, n-k + q)$-current in the usual sense on $\Omega$ that vanishes on forms in $\mathcal{N}_{k-p,k-q,Z}(\Omega)$.
    Conversely, if $T$ is a $(n-k + p, n-k + q)$-current on $\Omega$, that vanishes on forms in $\mathcal{N}_{k-p,k-q,Z}(\Omega)$,
    then $T$ defines a unique $(p,q)$-current $T'$ on $Z$ such that $i_* T' = T$.

    It is easy to see that the definitions of smooth forms, test forms and currents are independent of the embedding,
    and hence by gluing together in the same way one does on a complex manifold,
    we can define the sheaves of smooth forms, test forms and currents on any analytic space $Z$.
    Note in particular that by a smooth function on $Z$, we mean a function which is locally the restriction of a smooth
    function in the ambient space.

    In $\C$, one can define the principal value current $1/z^n = |z|^{2\lambda}/z^n |_{\lambda = 0}$ by analytic continuation,
    where $|_{\lambda = 0}$ denotes that for $\Re \lambda \gg 0$, we take the action of $|z|^{2\lambda}/z^n$ on a test form
    and take the value of the analytic continuation to $\lambda = 0$,
    which is easily seen to exist by a Taylor expansion, or integration by parts.
    Thus, if $\alpha$ is a smooth form on $\Cn$ and $\{ i_1,\cdots,i_m \} \subseteq \{ 1,\cdots,n \}$, with $i_j$ disjoint,
    then one gets a well-defined current
    \begin{equation} \label{elcurr}
        \frac{1}{z_{i_1}^{n_1}}\cdots\frac{1}{z_{i_k}^{n_k}}\dbar\frac{1}{z_{i_{k+1}}^{n_{k+1}}}\wedge\cdots\wedge\dbar\frac{1}{z_{i_m}^{n_m}}\wedge \alpha
    \end{equation}
    on $\Cn$ by taking $\dbar$ in the current sense together with tensor product of currents and multiplication of currents with smooth forms.
    In \cite{AW2}, if $\alpha$ has compact support, a current of the form \eqref{elcurr} is called an \emph{elementary current}.
    The class of pseudomeromorphic
    currents on a complex manifold was then introduced as currents that can be written as a locally finite sum of push-forwards of
    elementary currents. We will use the same definition on an analytic space $Z$.

    \begin{df}
        A current $\mu$ on $Z$ is said to be \emph{pseudomeromorphic}, denoted $\mu \in \PM(Z)$, if $\mu$
        can be written as a locally finite sum
        \begin{equation*}
            \mu = \sum (\pi_\alpha)_* \tau_\alpha,
        \end{equation*}
        where
        $\pi_\alpha : Z_\alpha \to Z$ is a family of compositions of modifications
        and open inclusions, and $\tau_\alpha$ are elementary currents on $Z_\alpha$.
    \end{df}

    Note in particular that, if $\pi : \tilde{Z} \to Z$ is a resolution of singularities of $Z$, and if $\mu \in \PM(\tilde{Z})$,
    then $\pi_* \mu \in \PM(Z)$.
    All the currents introduced in this article are pseudomeromorphic, as we will see directly from the proofs that the currents exist.
    In \cite{AW2}, it is shown that if $f$ is holomorphic on a complex manifold $X$, and $T \in \PM(X)$,
    one can define a multiplication $(1/f) T$ and $\dbar(1/f)\wedge T$. The same idea works equally well for strongly holomorphic
    functions on an analytic space.

    \begin{prop} \label{pmexist}
    Let $f$ be strongly holomorphic on
    $Z$, such that $f$ does not vanish on any irreducible component of $Z$, and let
    $T \in \PM(Z)$. Then the currents
    \begin{align*}
        \frac{1}{f} T := \left.\frac{|f|^{2\lambda}}{f}T \right|_{\lambda = 0} \quad \text{and} \quad \dbar\frac{1}{f}\wedge T :=
        \left.\frac{\dbar |f|^{2\lambda}}{f} \wedge T\right|_{\lambda = 0},
    \end{align*}
    where the right-hand sides are defined originally for $\Re \lambda \gg 0$, have current-valued analytic continuations to
    $\Re \lambda > -\epsilon$ for some $\epsilon > 0$, and the values at $\lambda = 0$ are pseudomeromorphic.
    The currents satisfies the Leibniz rule
    \begin{equation*}
        \dbar\left( \frac{1}{f} T\right) = \dbar \frac{1}{f} \wedge T + \frac{1}{f}\dbar T,
    \end{equation*}
    and $\supp (\dbar(1/f)\wedge T) \subseteq Z_f \cap \supp T$.
    If $f \neq 0$, then $(1/f) T$ defined in this way coincides with the usual multiplication of $T$ with the smooth function $1/f$.
    \end{prop}

    \begin{proof}
        If $Z$ is smooth, this is Proposition 2.1 in \cite{AW2}, except for the last statement.
        However, if $f \neq 0$, then $|f(z)|^{2\lambda}/f(z)$ is smooth
        in both $\lambda$ and $z$, and analytic in $\lambda$, so if $\xi$
        is a test form, $T . ((|f|^{2\lambda}/f)\xi)$ is analytic in $\lambda$, and hence the analytic continuation to $\lambda = 0$
        coincides with the value $T . ((1/f)\xi)$ at $\lambda = 0$.
        The proof in the general case goes through word for word as in the smooth case in Proposition 2.1 in \cite{AW2}.
    \end{proof}

    The crucial point in the proof of the following proposition is that for any analytic
    subset $W \subseteq Z$ and any $T \in \PM(Z)$, there exist natural restrictions
    \begin{equation} \label{eqrestroper}
     {\bf 1}_{W^c} T := |h|^{2\lambda} T |_{\lambda = 0} \quad \text{and} \quad {\bf 1}_W T := T - {\bf 1}_{W^c} T
    \end{equation}
    where $h$ is a tuple of holomorphic functions such that $W = \{ h = 0 \}$.
    The restrictions are independent of the choice of such $h$, and are such that $\supp {\bf 1}_{W} T \subseteq W$.
    This is Proposition 2.2 in \cite{AW2}, and the proof will go through in exactly the same way when $Z$ is an analytic space.

    \begin{prop} \label{pseudom0}
        Assume that $\mu \in \PM(Z)$, and that $\mu$ has support on a variety $V$. If $I_V$ is the ideal of holomorphic functions
        vanishing on $V$, then $\bar{I}_V \mu = 0$. If $\mu$ is of bidegree $(*,p)$, and $V$ has codimension $\geq p+1$ in $Z$, then $\mu = 0$.
    \end{prop}

    In the case that $Z$ is a complex manifold, this is Proposition 2.3 and Corollary 2.4 in \cite{AW2},
    and the proof there will go through in the same way also when $Z$ is an analytic space.
    The final step in the proof that $\mu = 0$ in the smooth case is to prove that $\mu = 0$ on $V_\reg$,
    which is proved with the help of the previous part of the proposition, and by degree reasons, and then by
    induction over the dimension of $V$, $\mu = 0$.
    In the singular case, this is done in the same way. Since this is a local statement, we can assume that $Z \subseteq \Omega \subseteq \Cn$,
    and consider $V$ as a subvariety of $\Omega$. Then, for the same reasons as in the smooth case, we get that $i_* \mu = 0$ on $V_\reg$,
    and by induction over the dimension of $V$ that $i_* \mu = 0$, and hence $\mu = 0$.

\section{Coleff-Herrera products of weakly holomorphic functions} \label{sectchp}

    Let $f_1,\cdots,f_{q+p} \in \Ow(Z)$.
    We want to define the Coleff-Herrera product
    \begin{equation*}
        T = \frac{1}{f_1}\cdots\frac{1}{f_q}\dbar\frac{1}{f_{q+1}}\wedge\cdots\wedge\dbar\frac{1}{f_{q+p}}.
    \end{equation*}
    If $f$ is strongly holomorphic, one way to define it is by
    \begin{equation} \label{tprim}
        T = \left.\frac{|f_1|^{2\lambda_1}\cdots|f_q|^{2\lambda_q}}{f_1\cdots f_q}
        \frac{\dbar |f_{q+1}|^{2\lambda_{q+1}}\wedge\cdots\wedge\dbar |f_{q+p}|^{2\lambda_{q+p}}}{f_{q+1}\cdots f_{q+p}}
        \right|_{\lambda_{q+p} = 0,\dots,\lambda_1 = 0},
    \end{equation}
    which a priori is defined only when $\Re \lambda_i \gg 0$; however, by Proposition~\ref{pmexist} it has an analytic continuation in
    $\lambda_{q+p}$ to $\Re \lambda_{q+p} > -\epsilon$
    for some $\epsilon > 0$, and the value at
    $\lambda_{q+p} = 0$
    is pseudomeromorphic.
    Again, by Proposition~\ref{pmexist}, it has an analytic continuation in
    $\lambda_{q+p-1}$ to $\lambda_{q+p-1} = 0$
    and so on, and hence the value at
    $\lambda_{q+p} = 0, \cdots, \lambda_1 = 0$
    exists.

    Note that if $\pi : Y \to Z$ is any modification of $Z$, we can define the corresponding Coleff-Herrera product of $f' = \pi^* f$ in
    $Y$. Taking the push-forward of this current to $Z$ will in fact give the Coleff-Herrera product of
    $f$ on $Z$. To see this, let $T^\lambda$ denote the form on the right-hand side of \eqref{tprim},
    with $\Re \lambda_i \gg 0$ fixed, and let ${T'}^{\lambda}$ denote the corresponding form on $Y$
    with $f'$ instead of $f$. If $\Re \lambda_i \gg 0$, then $T^\lambda$ and ${T'}^\lambda$ are smooth,
    and $\pi^* T^\lambda = {T'}^{\lambda}$, so $\pi_* {T'}^\lambda = T^\lambda$, since $\pi$ is a modification.
    Thus, by analytic continuation,
    $T = T^\lambda |_{\lambda_{q+p} = 0,\cdots,\lambda_1 = 0} = \pi_* {T'}^{\lambda} |_{\lambda_{q+p}=0,\cdots,\lambda_1 = 0} = \pi_* T'$.

    Now, if $f$ is weakly holomorphic, let $\pi : Z' \to Z$ be the normalization of $Z$, and $f' = \pi^* f$
    which is strongly holomorphic on $Z'$. Hence, the current
    \begin{equation} \label{eqtprim}
        T' = \frac{1}{f_1'}\cdots\frac{1}{f_q'}\dbar\frac{1}{f_{q+1}'}\wedge\cdots\wedge\dbar\frac{1}{f_{q+p}'}
    \end{equation}
    exists.
    \begin{df} 
        If $f = (f_1,\dots,f_{q+p})$
        is weakly holomorphic, we define the \emph{Coleff-Herrera product}
        \begin{equation} \label{chpdef}
            T = \frac{1}{f_1}\cdots\frac{1}{f_q}\dbar\frac{1}{f_{q+1}}\wedge\cdots\wedge\dbar\frac{1}{f_{q+p}}
        \end{equation}
        of $f$ as $\pi_* T'$,
        where $T'$ is defined by \eqref{eqtprim}.
    \end{df}
    If $f$ is strongly holomorphic, this definition will be the same as the definition in \eqref{tprim} since
    by the remark above, $T$ can be defined as the push-forward from any modification. In addition, if $f$ is weakly holomorphic,
    it can be defined by the push-forward of the corresponding current in any normal modification, since any normal
    modification factors through the normalization.

    We will call the factors $1/f_i$ the \emph{principal value factors}, and $\dbar (1/f_i)$ the \emph{residue factors}.
    \begin{remark} \label{remmixpvresfactors}
    Note that even though here, the principal value factors are to the left of the residue factors, 
    we could equally well have the residue and principal value factors mixed.
    However, changing the order will in general give a different current, but as we will see in Theorem~\ref{commfimult},
    if $f_i$ define a complete intersection, the current will not depend on the order (up to change of signs).
    \end{remark}

    \begin{remark} \label{epsilonlambda}
        The Coleff-Herrera product for $f = (f_1,\dots,f_p)$ strongly holomorphic is originally defined in \cite{CH}
        as the limit of integrals over $\cap \{ |f_i| = \epsilon_i(\delta) \}$ as $\epsilon \to 0$, where $\epsilon(\delta)$ tends to $0$ along
        an admissible path, cf., \eqref{chporig}.
        When $\epsilon(\delta)$ tends to $0$ along an admissible path, this will correspond to taking the analytic continuation to $\lambda = 0$
        in the order as in \eqref{tprim}, and in fact, for arbitrary $f$, the definition in \eqref{chporig} is equal to the one in \eqref{chpdef}
        defined by analytic continuation, see \cite{LS}.
    \end{remark}

    In \cite{Den} Denkowski gave a definition of the Coleff-Herrera product of $f$, for $f$ c-holomorphic,
    and we will see below that his definition coincides with ours in that case. The idea in \cite{Den} was to consider the graph of $f$,
    \begin{equation*}
        \Gamma_f = \{ (z,f(z)) \in Z\times \C^p_w | z \in Z \},
    \end{equation*}
    and even though $f$ is only c-holomorphic, the graph will be analytic.
    If $(z,w) \in \Gamma_f$, then $w = f(z)$, and hence on the graph $f_i = w_i$ is a strongly holomorphic function.
    If $\Pi$ is the projection from the graph to $Z$, since $f$ is continuous, $\Pi$ is a homeomorphism and in particular proper.
    The Coleff-Herrera product of $f$ was then defined by
    \begin{equation} \label{dendef}
        \dbar \frac{1}{f_1}\wedge \dots \wedge \dbar \frac{1}{f_p} = \Pi_*\left( \dbar\frac{1}{w_1}\wedge \dots \wedge \dbar \frac{1}{w_p} \right),
    \end{equation}
    and since $f_i = w_i$ on $\Gamma_f$, this should be a reasonable definition of the Coleff-Herrera product of $f$.
    The next proposition shows, as one might hope, that the definition of Denkowski coincides with ours.

    \begin{prop} \label{denequal}
        If $f = (f_1,\dots,f_p)$ is c-holomorphic, then the definition of the Coleff-Herrera product of $f$ in \eqref{chpdef}
        and in \eqref{dendef} coincide.
    \end{prop}

    \begin{proof}
        In \cite{Den} the definition used for the Coleff-Herrera product of strongly holomorphic functions
        was the one from \cite{CH}. However, by Remark~\ref{epsilonlambda}
        we can assume that the definition by analytic continuation is used instead.
        Let $\pi : Z' \to Z$ be the normalization of $Z$ and $f' = \pi^* f$. We have projections
        $\Pi : \Gamma_f \to Z$ and $\Pi' : \Gamma_{f'} \to Z'$, where $\Gamma_f \subseteq Z\times \C_w^p$ and
        $\Gamma_{f'} \subseteq Z' \times \C_{w'}^p$ are the graphs of $f$ and $f'$.
        Thus we have a commutative diagram
        \begin{equation} \label{moddiag}
        \begin{CD}
            \Gamma_{f'} @>(\pi\times \Id)|_{\Gamma_{f'}} >> \Gamma_f \\
            @VV \Pi' V @VV \Pi V \\
            Z' @> \pi >> Z.
        \end{CD}
        \end{equation}
        We will denote the current
        $\dbar (1/f_1')\wedge \dots \wedge \dbar (1/f_p')$
        on $Z'$ by $\mu^{f'}$, and similarly for $\mu^{w}$
        and $\mu^{w'}$ defined on $\Gamma_f$ and $\Gamma_{f'}$ respectively.
        Then
        $\dbar (1/f_1)\wedge\dots\wedge \dbar (1/f_p)$
        is defined in \eqref{chpdef} by $\pi_* \mu^{f'}$,
        and in \eqref{dendef} by $\Pi_* \mu^{w}$.
        Now, $(\pi\times \Id)|_{\Gamma_{f'}} : \Gamma_{f'} \to \Gamma_f$ is a modification of $\Gamma_f$ so we have
        $\mu^w = (\Pi\times\Id)_* \mu^{w'}$,
        and since $\Pi' : \Gamma_{f'} \to Z'$ is a biholomorphism and $w_i' = f_i'$ on $\Gamma_{f'}$ we also have
        $\mu^{f'} = \Pi'_* \mu^{w'}$.
        Thus both are the push-forward of the same current in $\Gamma_{f'}$, and since the diagram \eqref{moddiag} commutes,
        both will have the same push-forward to $Z$.
    \end{proof}

    The next two theorems are extensions to the case of weakly holomorphic functions of well-known results
    of the Coleff-Herrera product of strongly holomorphic functions (in the case $q = 0$ or $q = 1$),
    see \cite{CH}, or the case of holomorphic functions on a complex manifold, see \cite{P1}.

    \begin{thm} \label{leibnizsupport}
        If
        $f=(f_1,\cdots,f_{q+p})$
        is weakly holomorphic, then $T$, defined by \eqref{chpdef}, satisfies the Leibniz rule
    \begin{equation*}
        \dbar T = \sum_{j= 1}^q
        \frac{1}{f_1}\cdots\frac{1}{f_{j-1}}\dbar\frac{1}{f_j}\wedge\frac{1}{f_{j+1}}\cdots\frac{1}{f_q}
        \dbar\frac{1}{f_{q+1}}\wedge\cdots\wedge\dbar\frac{1}{f_{q+p}},
    \end{equation*}
        and $\supp T \subseteq Z_{(f_{q+1},\cdots,f_{q+p})}$.
    \end{thm}

\begin{proof}
    First we assume that $f$ is strongly holomorphic. Then
    the Leibniz rule follows by analytic continuation, since if $\Re \lambda \gg 0$, we have
\begin{equation*}
        \dbar\left(\frac{|f|^{2\lambda}}{f}\right) = \frac{\dbar |f|^{2\lambda}}{f} \text{ and }
    \dbar\left(\frac{\dbar |f|^{2\lambda}}{f}\right) = 0.
\end{equation*}
    The weakly holomorphic case follows by taking push-forward from the normalization.
    For the last part, let $T'$ be the current corresponding to $T$ in the normalization, and $f' = \pi^* f$ be the pull-back of
    $f$ to the normalization. Then by Proposition~\ref{pmexist},
    $T' = 0$ outside of $Z_{f_i'}$,
    $i \geq q+1$,
    and hence
    $\supp T \subseteq \pi(\supp T' ) \subseteq \pi(Z_{(f_{q+1}',\cdots,f_{q+p}')}) = Z_{(f_{q+1},\cdots,f_{q+p})}$,
    where the last equality follows from Proposition~\ref{zeroproperties}.
\end{proof}

    It is natural in this context to ask how to define a reasonable multiplication of a weakly holomorphic function
    with a current, something which we will need in the case that the current is a Coleff-Herrera product
    to be able to state the next theorem.
    If $g \in \Ow(Z)$, and $T$ is the Coleff-Herrera product in \eqref{chpdef}, we define $g T$ by
    \begin{equation} \label{wholomult}
        g T = \pi_*(\pi^*g T'),
    \end{equation}
    where $\pi : Z' \to Z$ is the normalization of $Z$, and $T'$ is the corresponding Coleff-Herrera product of $f' = \pi^* f$.
    In the case that both $f$ and $g$ are c-holomorphic, Denkowski gives a definition of multiplication
    of $g$ and the Coleff-Herrera product of $f$ in \cite{Den},
    and by a similar argument as that in Proposition~\ref{denequal}, one sees that our definition
    coincides with the one in \cite{Den} in that case.
    Note however, that we do not define a multiplication of a weakly holomorphic function
    with an arbitrary current, and as we will see in Section~\ref{sectmult},
    this will not be possible if we require it to satisfy certain natural properties.

    \begin{thm} \label{commfimult}
    Let
    $f = (f_1,\dots,f_{q+p})$
    be weakly holomorphic, such that
    $(f_{q+1},\dots,f_{q+p})$
    defines a complete intersection, and that
    $(f_i,f_{q+1},\dots,f_{q+p})$
    defines a complete intersection for
    $1 \leq i \leq q$.
    Then the principal value factors in
    \begin{equation*}
        T = \frac{1}{f_1}\cdots\frac{1}{f_q}\dbar\frac{1}{f_{q+1}}\wedge\cdots\wedge\dbar\frac{1}{f_{q+p}}
    \end{equation*}
    commute with other principal value factors or residue
    factors (see Remark~\ref{remmixpvresfactors}),
    and the residue factors anticommute.
    In addition, if
    $1 \leq k \leq q$,
    we have
    \begin{equation} \label{fimulteq1}
        f_k T = \frac{1}{f_1}\cdots\widehat{\frac{1}{f_k}}\cdots\frac{1}{f_q}\dbar\frac{1}{f_{q+1}}\wedge\cdots\wedge\dbar\frac{1}{f_{q+p}},
    \end{equation}
    and if
    $q+1 \leq j \leq q+p$,
    then
    \begin{equation} \label{fimulteq2}
        f_j T = 0.
    \end{equation}
    \end{thm}

Note that in case $f_i \in \Ow(Z)$, then the left-hand sides of \eqref{fimulteq1} and \eqref{fimulteq2}
are defined by \eqref{wholomult}.

    \begin{remark}
        In the smooth case, the first part of Theorem~\ref{commfimult} (about permuting the factors)
        follows from the theorem of Samuelsson in \cite{S}, about the analyticity of the residue integral \eqref{tprim}.
        In fact, his theorem holds also for strongly holomorphic functions on an analytic space, cf., \cite{LS}.
        Since the proof of the first part of Theorem~\ref{commfimult} reduces to the strongly holomorphic
        case, one could thus refer to the results of Samuelsson. However, since the proof of this deep theorem
        of Samuelsson is quite involved, we still prefer to give a direct proof of the first part of Theorem~\ref{commfimult},
        since it can be done by much more elementary means.
    \end{remark}

Note that in the following lemmas, which we will use to prove Theorem~\ref{commfimult},
we assume that the functions are strongly holomorphic.

\begin{lma} \label{comm}
    Assume that $f_1, f_2 \in \Ok(Z)$ and that $T \in \PM(Z)$ is of bidegree $(*,p)$.
    If $Z_{f_1}\cap Z_{f_2} \cap \supp T \subseteq V$, for some analytic set $V \subseteq Z$ of codimension $\geq p+1$ in $Z$,
    then
    \begin{equation} \label{commeq1}
        \frac{1}{f_1}\frac{1}{f_2} T = \frac{1}{f_2}\frac{1}{f_1} T.
    \end{equation}
    If $Z_{f_1}\cap Z_{f_2}\cap \supp T \subseteq V'$, for some analytic set $V'$ of codimension $\geq p+2$ in $Z$,
    then 
    \begin{equation} \label{commeq2}
        \frac{1}{f_1} \dbar\frac{1}{f_2} \wedge T = \dbar\frac{1}{f_2}\wedge \frac{1}{f_1}T,
    \end{equation}
    and if in addition $Z_{f_1} \cap Z_{f_2}\cap \supp \dbar T \subseteq V''$, for some analytic set $V''$ of codimension $\geq p+3$, then
    \begin{equation} \label{commeq3}
        \dbar \frac{1}{f_1}\wedge\dbar \frac{1}{f_2} \wedge T = - \dbar \frac{1}{f_2}\wedge\dbar\frac{1}{f_1}\wedge T.
    \end{equation}
\end{lma}

\begin{proof}
    We have by Proposition~\ref{pmexist} that
    \begin{equation} \label{eqf1f2T}
        \frac{1}{f_1}\frac{1}{f_2} T - \frac{1}{f_2}\frac{1}{f_1} T,
    \end{equation}
    is zero outside of $Z_{f_1}$, since both terms
    are just multiplication of $(1/f_2)T$ with the smooth function $(1/f_1)$, and similarly
    it is zero outside of $Z_{f_2}$. Thus \eqref{eqf1f2T}
    is a pseudomeromorphic current on $Z$ of bidegree $(*,p)$ with support on $Z_{f_1} \cap Z_{f_2} \cap V$, which has codimension $\geq p+1$, so \eqref{commeq1} follows by Proposition~\ref{pseudom0}.
    Similarly outside of $Z_{f_1}$, we get
    that
    \begin{equation} \label{eqf1dbarf2T}
        \frac{1}{f_1} \dbar \frac{1}{f_2} \wedge T - \dbar \frac{1}{f_2}\wedge \frac{1}{f_1} T
    \end{equation}
    is zero, so \eqref{eqf1dbarf2T}
    is a pseudomeromorphic current on $Z$ of bidegree $(*,p+1)$ and has support on $Z_{f_1} \cap Z_{f_2} \cap \supp T$,
    so \eqref{commeq2} follows by Proposition~\ref{pseudom0}.
    For \eqref{commeq3}, we get by Theorem~\ref{leibnizsupport} and \eqref{commeq2} that
    \begin{align*}
        &\dbar \frac{1}{f_1} \wedge \dbar \frac{1}{f_2} \wedge T = \dbar\left( \frac{1}{f_1} \dbar \frac{1}{f_2} \wedge T\right) +
        \frac{1}{f_1} \dbar \frac{1}{f_2} \wedge\dbar T \\
        &=  \dbar\left(\dbar \frac{1}{f_2} \wedge \frac{1}{f_1} T\right) +  \frac{1}{f_1} \dbar \frac{1}{f_2} \wedge\dbar T \\
        &= -\dbar \frac{1}{f_2} \wedge\dbar \frac{1}{f_1} \wedge T - \dbar \frac{1}{f_2} \wedge \frac{1}{f_1} \dbar T +
        \frac{1}{f_1} \dbar \frac{1}{f_2} \wedge\dbar T 
        = -\dbar \frac{1}{f_2} \wedge\dbar \frac{1}{f_1} \wedge T
    \end{align*}
    where the last equality holds because of \eqref{commeq2} and the assumption of the support of $\dbar T$.
\end{proof}

\begin{lma} \label{mult2}
    Assume $f, g \in \Ok(Z)$, and $f/g \in \Ok(Z)$. If $T \in \PM(Z)$ has bidegree $(*,p)$ and 
    $Z_g\cap \supp T \subseteq V$, for some analytic subset $V$ of codimension $\geq p+1$, then
\begin{equation*}
    f \left(\frac{1}{g} T\right) = \frac{f}{g} T.
\end{equation*}
\end{lma}

\begin{proof}
    Outside of $Z_g$, we can see $(1/g) T$ as multiplication by the smooth function $1/g$
    by Proposition~\ref{pmexist}. Hence we have $f (1/g) T = (f/g) T$ since their difference is a pseudomeromorphic
    current with support on $Z_g \cap \supp T$, so it is $0$ by Proposition~\ref{pseudom0}.
\end{proof}

\begin{proof} [Proof of Theorem~\ref{commfimult}]
    First we observe that it is enough to prove the theorem in case $f_i$ are strongly holomorphic, since if $\pi : Z' \to Z$
    is the normalization of $Z$, and $f' = \pi^* f$, then $f'$ is a complete intersection, and if the theorem
    holds in $Z'$, it holds in $Z$ by taking push-forward of the corresponding currents.
    Hence, we can assume that $f_i \in \Ok(Z)$, and the commutativity properties will then follow from Lemma~\ref{comm}.
    For example, if we want to see that
    $1/f_i$ and $1/f_{i+1}$
    commute, we can apply Lemma~\ref{comm} with
    \begin{equation*}
        T = \frac{1}{f_{i+2}}\cdots\frac{1}{f_q}\dbar\frac{1}{f_{q+1}}\wedge\cdots\wedge\dbar\frac{1}{f_{q+p}},
    \end{equation*}
    and then multiply with
    $(1/f_1)\cdots(1/f_{i-1})$
    from the left. In case some of the residue factors, say
    $f_{q+1},\ldots,f_{q+k}$,
    are to the left of the principal value factors, then
    $Z_{(f_{q+k+1},\ldots,f_{q+p})}$ has codimension $p-k$
    in a neighborhood of $Z_f \supseteq \supp T$ and the result follows in the same way from Lemma~\ref{comm}.
    The other cases follow similarly from Lemma~\ref{comm}.

    The equality \eqref{fimulteq1} follows from Lemma~\ref{mult2} since $Z_f$
    has codimension $p$. By the first part of the theorem, we can assume that
    $j = q+1$
    in \eqref{fimulteq2}. Then
    \begin{align*}
        f_{q+1}\left(\dbar\frac{1}{f_{q+1}}\wedge\cdots\wedge\dbar\frac{1}{f_{q+p}}\right) &= 
        \dbar\left(f_{q+1}\frac{1}{f_{q+1}}\wedge\dbar\frac{1}{f_{q+2}}\wedge\cdots\wedge\dbar\frac{1}{f_{q+p}}\right) \\
        &= \dbar\left(\dbar\frac{1}{f_{{q+2}}}\wedge\cdots\wedge\dbar\frac{1}{f_{q+p}}\right) = 0
    \end{align*}
    by \eqref{fimulteq1}, and Theorem~\ref{leibnizsupport}.
\end{proof}

\section{Multiplication of currents with weakly holomorphic functions} \label{sectmult}

Now, we will return to the issue of multiplication of currents with weakly holomorphic functions. Assume $g \in \Ow(Z)$
and $S \in \PM(Z)$. Since $S \in \PM(Z)$, we have $S = \sum (\pi_\alpha)_* \tau_\alpha$, where $\tau_\alpha$ are elementary currents
on the complex manifolds $Z_{\alpha}$. Given such a decomposition, since any normal modification of $Z$ factors through the
normalization, that is, $\pi_\alpha = \pi \circ \nu_\alpha$, for some $\nu_\alpha : Z_\alpha \to Z'$,
we get a current $S'$ in the normalization $Z'$ of $Z$ such that $\pi_* S' = S$ by taking the
push-forward of $\tau_\alpha$ to $Z'$, i.e., $S' = \sum (\nu_\alpha)_* \tau_\alpha$.
To define multiplication of the Coleff-Herrera product with
the weakly holomorphic function $g$ in \eqref{wholomult}, we defined it as the push-forward of $\pi^* g S'$.
In general, the current $S'$ will depend on the decomposition $S = \sum (\pi_\alpha)_* \tau_\alpha$.
However, in \eqref{wholomult} we had a canonical representative in the normalization,
and hence the multiplication was well-defined.
The following example however shows that this multiplication depends on this choice of representative.

\begin{ex} \label{exmultcholo}
Let $\pi : \C^n \to \C^{2n}$ be defined by
\begin{equation*}
    \pi(t_1,\cdots,t_n) = (t_1,\cdots,t_{n-1},t_1^2t_n,\cdots,t_{n-1}^2t_n,t_n^2,t_n^5).
\end{equation*}
Then $\pi$ is proper and injective, so $\pi(\C^n) = Z$ is an analytic variety of dimension $n$. Since $(\partial\pi_j/\partial z_i)_{i,j}$
has full rank outside of $\{ 0 \}$, $Z_\sing \subseteq \{ 0 \}$, and we will see below that actually $Z_\sing = \{ 0 \}$. Let
\begin{equation*}
    \tilde{S} = \dbar\frac{1}{t_1}\wedge\cdots\wedge\dbar\frac{1}{t_{n-1}}\wedge\dbar\frac{1}{t_n^3}
\end{equation*}
and $S = \pi_* \tilde{S}$.
Then, since $d(t_nt_i^2) = t_i(2t_ndt_i + t_idt_n)$ and $dt_n^5 = 5t_n^4dt_n$, $dz_k \wedge S = 0$ for $k = n,\dots,2(n-1)$, and $dz_{2n} \wedge S = 0$.
Hence if $S . \xi \neq 0$, then $\xi$ must be of the form $\xi = \xi_0 dz_1\wedge\cdots\wedge dz_{n-1}\wedge dz_{2n-1}$. We have
\begin{align*}
    &S . \xi = \tilde{S}.\xi_0 dt_1\wedge\cdots\wedge dt_{n-1}\wedge 2t_n dt_n = \\
    &2\cdot (2\pi i)^n \left.\left(\sum_{i=1}^{n-1} t_i^2\frac{\partial}{\partial z_{n-1 + i}}\xi_0 + 2t_n\frac{\partial}{\partial z_{2n-1}}\xi_0 +
    5t_n^4\frac{\partial}{\partial z_{2n}}\xi_0 \right)\right|_{t = 0} = 0,
\end{align*}
and thus $S = 0$.
However,
\begin{equation*}
    t_n \tilde{S}.\xi dt_1\wedge\cdots\wedge dt_n^2 = 2(2\pi i)^n\xi(0)
\end{equation*}
so $\pi_*(t_n \tilde{S}) = \pi_*( \pi^* g \tilde{S}) \neq 0$, where $g \in \Oc(Z)$ is such that $\pi^*g = t_n$.
Note that $g$ is not strongly holomorphic at $0$, and hence $Z_\sing = \{ 0 \}$.
Thus, since $S = \pi_* \tilde{S} = 0$, while $\pi_*( (\pi^* g) \tilde{S}) = 0$,
it is impossible to define a multiplication of currents with weakly holomorphic functions
in a way compatible with push-forwards, i.e., that $gS$ only depends on $g$ and $S$, and
such that $g S = \pi_*( (\pi^*g) S')$ if $S = \pi_* S'$.
\end{ex}

Hence, the multiplication in \eqref{wholomult} does not depend only on $g$ and $S$, but also on the functions
$f$ defining $S$.
Recall that the \emph{pole set}, $P_\phi$, of a meromorphic function $\phi$ is the
set where $\phi$ is not strongly holomorphic.
Recall also the definitions of the restriction operators ${\bf 1}_V$ and ${\bf 1}_{V^c}$ in \eqref{eqrestroper}.
If we require that the current we get in the multiplication has restriction $0$ to $P_\phi$, the multiplication
is in fact uniquely defined in $\PM(Z)$, as the following proposition shows.
This can in some cases be a natural condition, and in fact even automatic in some cases, see Corollary~\ref{corwholomult3}.
However, in Example~\ref{exmultcholo}, since the common zero set of the functions defining $S$
equals the pole set of $g$, we expect $S$ and $gS$ to have its support on $P_g$, and hence
the condition is not very natural then.

\begin{prop} \label{wholomult2}
    Let $\mu \in \PM(Z)$ and $\phi \in \Ow(Z)$.
    Then, there exists a unique current, denoted $\phi \mu$, in $\PM(Z)$, such that $\phi \mu$ is
    just multiplication of the smooth function $\phi$ with the current $\mu$ outside of $P_\phi$,
    and ${\bf 1}_{P_\phi}(\phi \mu ) = 0$. If $\mu = \pi_* \mu'$, where $\pi : Z' \to Z$ is the normalization of
    $Z$ and $\mu' \in \PM(Z')$, then
    \begin{equation} \label{eqphimudef}
        \phi \mu = \pi_* ( (\pi^* \phi) {\bf 1}_{(\pi^{-1}(P_\phi))^c}\mu').
    \end{equation}
\end{prop}

\begin{proof}
   First, we prove the uniqueness. Assume that $T_1$ and $T_2$ are two such currents, so that
   $T_1 - T_2$ has support on $P_\phi$. Hence, ${\bf 1}_{P_\phi^c} (T_1 - T_2) = 0$.
   But then,
   \begin{equation*}
      T_1 - T_2 = {\bf 1}_{P_\phi^c}(T_1 - T_2) + {\bf 1}_{P_\phi}(T_1 - T_2) = 0,
   \end{equation*}
   since ${\bf 1}_{P_\phi} T_1 = {\bf 1}_{P_\phi} T_2 = 0$.
   Thus, we only need to prove that $\phi \mu$ in \eqref{eqphimudef} satisfies the conditions
   in the proposition. It is clear that the right-hand side of \eqref{eqphimudef} is just multiplication
   of $\phi$ with $\mu$ outside of $P_\phi$. Hence, it remains to prove that ${\bf 1}_{P_\phi}(\phi\mu) = 0$.
   However,
   \begin{equation*}
      {\bf 1}_{P_\phi}(\phi\mu) = \pi_*( {\bf 1}_{\pi^{-1}(P_\phi)}(\pi^*\phi) {\bf 1}_{(\pi^{-1}(P_\phi))^c} \mu') = 0,
   \end{equation*}
   since ${\bf 1}_V {\bf 1}_{V^c} = {\bf 1}_V(1-{\bf 1}_V)= 0$ because ${\bf 1}_V {\bf 1}_V = {\bf 1}_V$,
   and ${\bf 1}_V$ commutes with multiplication with smooth functions.
\end{proof}

\begin{cor} \label{corwholomult3}
    Assume that $\mu \in \PM(Z)$ is of bidegree $(*,p)$ and $\phi \in \Ow(Z)$ is such that
    $P_\phi$ has codimension $\geq p+1$ in $Z$. Then there exists a unique current $\phi \mu \in \PM(Z)$
    such that $\phi \mu$ coincides with the usual multiplication of $\phi$ with $\mu$ outside of $P_\phi$.
    If $\mu = \pi_* \mu'$, where $\pi : Z' \to Z$ is the normalization of $Z$ and $\mu' \in \PM(Z')$, then
    \begin{equation} \label{eqphimudef2}
        \phi \mu = \pi_*( (\pi^* \phi)\mu').
    \end{equation}
\end{cor}

\begin{proof}
    By Proposition~\ref{wholomult2}, the only thing we need to prove is that for any
    $T \in \PM(Z)$ and $T' \in \PM(Z')$ of bidegree $(*,p)$, we have ${\bf 1}_{P_\phi} T = 0$
    and ${\bf 1}_{\pi^{-1}(P_\phi)} T' = 0$.
    However, since $P_\phi$ has codimension $\geq p+1$, $\pi^{-1}(P_\phi)$ has codimension $\geq p+1$
    by Lemma~\ref{dimnormalization}. Hence, ${\bf 1}_{P_\phi} T = 0$ and ${\bf 1}_{\pi^{-1}(P_\phi)} T' = 0$
    by Proposition~\ref{pseudom0}, since the currents have support on $P_\phi$ and $\pi^{-1}(P_\phi)$
    respectively.
\end{proof}

Note, in particular that if $Z_\sing$ has codimension $\geq p+1$, the condition of the codimension of $P_g$ is automatically
satisfied for any weakly holomorphic function $g \in \Ow(Z)$.

Another question is whether the Coleff-Herrera product could be defined
as the analytic continuation of an integral on $Z$ rather than $Z'$.
A natural way to do this would be to try to regularize in \eqref{chpdef} by factors
$\dbar |F_i|^{2\lambda_i}$ instead of $\dbar |f_i|^{2\lambda_i}$,
where $F_i$ is a tuple of strongly holomorphic functions such that $Z_{F_i} = P_{1/f_i}$.
However, the analytic continuation to $\lambda = 0$ will in general
not coincide with our definition, even if $f$ defines a complete intersection, as the following example shows.
\begin{ex} \label{exchintrinsic}
    Let $Z = \{ z \in \C^3\ \lvert \ z_1^3 = z_2^2 \} = V\times \C$, which has normalization $\pi(s,t) = (s^2,s^3,t)$,
    and let
    $\pi^* f_1 = (1 + s)t$ and $\pi^* f_2 = s^2$.
    Then $Z_f = \{ 0 \}$, so $f$ is a complete intersection.
    Note that
    $\pi^*(1/f_1) = (1/t)(1-s + O(s^2))$
    for $|s|<1$, and that holomorphic functions in $s$ at the origin
    correspond to strongly holomorphic functions on $V$ at the origin precisely when the Taylor expansion at the origin
    contains no term $s$. Thus
    $P_{1/f_1} = \pi(\{ s = 0 \} \cup \{ s = -1\} \cup \{ t = 0 \})$,
    so if $\{ F = 0 \} \supseteq P_{1/f_1}$, then $\{ F = 0 \} \supseteq Z_{f_2}$.
    Thus $(\dbar |F|^{2\lambda} /f_1) \wedge \dbar (1/f_2) = 0$
    for $\Re \lambda \gg 0$. However, we have
    \begin{align*}
        &\dbar \frac{1}{f_1}\wedge\dbar\frac{1}{f_2} . \varphi dz_1\wedge dz_3 =
        \frac{1}{1 + s}\dbar\frac{1}{t}\wedge \dbar \frac{1}{s^2} . \varphi(s^2,s^3,t)ds^2\wedge dt = 4\pi i \varphi(0), \\
    \end{align*}
    so $\dbar (1/f_1)\wedge \dbar (1/f_2)$ is non-zero.
\end{ex}

\section{Bochner-Martinelli type residue currents} \label{sectbm}

We will show that we can define a Bochner-Martinelli type residue current associated with a tuple of weakly holomorphic functions,
either by using a similar approach as for the Coleff-Herrera product with the help of the normalization,
or by defining it intrinsically on $Z$ by means of analytic continuation. In view of Example~\ref{exchintrinsic},
it is not clear how to do this directly for the Coleff-Herrera product. In addition, we will show that for weakly holomorphic functions defining
a complete intersection, the Coleff-Herrera product and the Bochner-Martinelli current coincide, Theorem~\ref{bmch}.

Let $f = (f_1,\dots,f_p)$ be weakly holomorphic. 
We will follow the approach by Andersson from \cite{A1}, and make the identification $f = \sum f_i e_i^*$, where $(e_1,\cdots,e_p)$ is a frame
for a trivial vector bundle $E$ over $Z$, and $(e_1^*,\dots,e_p^*)$ is the dual frame.
Since we will only use the case of trivial vector bundles, this identification
is not strictly necessary. However, we use this since it greatly simplifies the notation in the proof of Lemma~\ref{uguf}.
Then, on the set where $f$ is strongly holomorphic, $\nabla_f := \delta_f - \dbar$ induces a complex on currents on $Z$ with
values in $\bigwedge E$, where $\delta_f$ is interior multiplication with $f$.
To construct the Bochner-Martinelli current we define
\begin{equation} \label{sigmau}
    \sigma = \sum \frac{\bar{f}_i e_i}{|f|^2} \quad \text{and} \quad u = \sum_{k=0}^{p-1} \sigma\wedge(\dbar\sigma)^k.
\end{equation}
Note that outside of $Z_f \cup P_{f_1}\cup\dots\cup P_{f_p}$, both $u$ and $\sigma$ are smooth, and $\nabla_f u = 1$.

Recall that a \emph{universal denominator} at a germ $(Z,z)$ is a strongly holomorphic function $h$, not
vanishing on any irreducible component of $(Z,z)$ such that $h \Ow_{Z,z} \subseteq \Ok_{Z,z}$.
For each $z \in Z$, there always exist a universal denominator $h$, such that $h$ is a universal denominator
in a neighborhood of $z$, see for example \cite{G}, Theorem~Q.2.

\begin{prop} \label{bmexist}
    Assume that $f = (f_1,\cdots,f_p)$ is weakly holomorphic on $Z$. Let $F$ be a tuple of strongly holomorphic functions,
    such that $\{ F = 0 \} \supseteq Z_f$, and $\{ F = 0 \}$ does not contain any irreducible component of $Z$,
    and let $h$ be a universal denominator on $Z$.
    Then the forms $|hF|^{2\lambda} u$ and $\dbar |hF|^{2\lambda}\wedge u$
    are arbitrarily smooth if $\Re \lambda \gg 0$, and have current-valued analytic continuations
    to $\Re \lambda > -\epsilon$ for some $\epsilon > 0$. The currents
    \begin{equation} \label{UR}
        U^f := |hF|^{2\lambda} u|_{\lambda = 0} \quad \text{and} \quad R^f := \dbar|hF|^{2\lambda} \wedge u|_{\lambda = 0}
    \end{equation}
    are independent of the choice of $F$ and $h$, and if $\pi : Y \to Z$ is a modification of $Z$, then $U^f = \pi_* U^{\pi^*f}$
    and $R^f = \pi_* R^{\pi^*f}$.
\end{prop}

\begin{proof}
    We first show that $|hF|^{2\lambda} u$ and $\dbar |hF|^{2\lambda} \wedge u$ are arbitrarily
    smooth when $\Re \lambda \gg 0$. Since $\dbar |hF|^{2\lambda} = |hF|^{2(\lambda-1)} \dbar |hF|^2$,
    it is enough to prove this for $|hF|^{2\lambda} u$. We let $g_i := h f_i$, where $g_i \in \Ok(Z)$ since
    $h$ is a universal denominator. If we differentiate $u$ outside of $\{ h = 0 \} \cup Z_f$, we get
    terms of the form $\xi/(h^k |f|^{2n})$, where $\xi$ is smooth, since if $f_i = g_i/h$, the terms
    in $u$ are smooth except for factors $h$ and $|f|^2$ in the denominators. Thus, we only need
    to see that $|hF|^{2\lambda}/(h^k |f|^{2n})$ tends to $0$ on $\{ h = 0 \} \cup Z_f$.
    This is clear outside of $Z_f$ if $\Re \lambda \gg 0$, so we need to prove that $|hF|^{2\lambda}/|f|^{2n}$
    tends to $0$ on $Z_f$.
    If we multiply the numerator and denominator by $|h|^{2n}$, we get
    \begin{equation} \label{equh}
        |h|^{2n}|hF|^{2\lambda}/(|hf|^{2n}).
    \end{equation}
    We note that $hf$ is strongly holomorphic, and in fact, $\{ hF = 0 \} \supseteq \{ hf = 0 \}$ because
    \begin{equation*}
        Z_{hf}  = \pi( Z_{\pi^*(hf)}) = \pi(Z_{\pi^* h}) \cup \pi(Z_{\pi^* f}) = Z_h \cup \pi(Z_{\pi^*f}) = \{ h = 0 \} \cup Z_f,
    \end{equation*}
    by Proposition~\ref{zeroproperties} and the fact that $\pi$ is surjective.
    Thus, \eqref{equh} will tend to $0$ on $Z_f$ by the Nullstellensatz
    if $\Re \lambda \gg 0$.

    Now, we assume that $Z$ is smooth. Then we can take $F = f$ and $h \equiv 1$, and in that case,
    the proposition is the existence part of Theorem 1.1 in \cite{A1}, except for the fact that
    $U^f = \pi_* U^{\pi^* f}$ and $R^f = \pi_* R^{\pi^* f}$, which however easily follows by analytic
    continuation. To see that the definition of $R^f$ is independent of the choice of $F$, we see from
    the proof of Theorem 1.1 in \cite{A1} that $\dbar |F|^{2\lambda} \wedge u$ acting on a test form $\varphi$
    becomes, with a suitable resolution of singularities $\pi : \tilde{X} \to X$,
    a finite sum of terms of the kind
    \begin{equation} \label{bmeqn}
        \int \frac{\dbar |u \mu_1|^{2\lambda}}{\mu_2} \wedge \sigma' \wedge \pi^*\varphi,
    \end{equation}
    where $\mu_1$ and $\mu_2$ are monomials such that $\{ \mu_1 = 0 \} \supseteq \{ \mu_2 = 0 \}$, $u$ is non-zero and $\sigma'$ is smooth.
    Thus, it is enough to observe that the value at $\lambda=0$ of \eqref{bmeqn} is independent of $\mu_1$
    (where $u\mu_1$ is the pull-back of $F$), as long as $\{ \mu_1 = 0 \} \supseteq \{ \mu_2 = 0 \}$.
    In the same way, one sees that the definition of $U^f$ is independent of the choice of $F$.

    Now, if $f$ is weakly holomorphic, and $\pi : \tilde{Z} \to Z$ is a resolution of singularities,
    from the smooth case we know that $\dbar |\pi^*(hF)|^{2\lambda} \wedge \pi^* u$ has a current-valued analytic
    continuation to $\lambda = 0$ independent of the choice of $hF$.
    Hence, the weakly holomorphic case follows by taking push-forward, since
    $\dbar |hF|^{2\lambda} \wedge u = \pi_*(\dbar |\pi^* (hF)|^{2\lambda} \wedge \pi^* u)$ for $\Re \lambda \gg 0$.
\end{proof}

In fact, to prove the existence of $U^f$ and $R^f$, defined by \eqref{UR}, it is sufficient to use $|F|^{2\lambda} u$
and $\dbar |F|^{2\lambda} \wedge u$, which can be seen are integrable on $Z$ if $\Re \lambda \gg 0$ by going back to the normalization.
However, the addition of the universal denominator $h$ ensures that the forms are (arbitrarily) smooth if $\Re \lambda \gg 0$.

The following properties of the Bochner-Martinelli current, $R^f$, are well-known in the smooth case, see \cite{PTY} and \cite{A1}.

\begin{prop} \label{bmprop}
    Let $f = (f_1,\cdots,f_p)$ be weakly holomorphic, and assume that $p' = \codim Z_f$.
    The current $R^f$ has support on $V = Z_f$, and 
    there is a decomposition $R^f = \sum_{k = p'}^p R_k$, where $R_k \in \PM(Z)$ is a $(0,k)$-current with values in $\bigwedge^k E$.
    In addition, if $f$ is strongly holomorphic, then $R^f = 1-\nabla_{f} U^f$.
\end{prop}

\begin{proof}
    In case $Z$ is a complex manifold, this is parts of Theorem~1.1 in \cite{A1}, except for the fact
    that $R_k \in \PM(Z)$. However, that $R_k$ is pseudomeromorphic can, as was noted in \cite{AW2}, easily be seen from the proof
    of Theorem 1.1 in \cite{A1}. The proposition then follows in case of an analytic space, by taking push-forward from
    a resolution of singularities, except for the fact that $R^f = \sum_{k= p'}^p R_k$, where $p' = \codim Z_f$,
    since modifications does not in general preserve codimensions of subvarieties. However, we get that $R^f = \sum_{k=0}^p R_k$, where
    $R_k \in \PM(Z)$ is a $(0,k)$-current, and $R_k$ has support on $Z_f$. Thus, by Proposition~\ref{pseudom0}, $R_k = 0$
    for $k < \codim Z_f = p'$.
\end{proof}

\begin{remark} \label{rembmintrinsic}
    If the mapping $f$ is weakly holomorphic, as we saw in Example~\ref{exmultcholo}, we do not have a well-defined
    multiplication of weakly holomorphic functions with pseudomeromorphic currents on $Z$.
    Hence, the formula $R^f = 1 - \nabla_f U^f$ in the strongly holomorphic case does not necessarily
    have any meaning if $f$ is weakly holomorphic.
    However, one can give this multiplication meaning by Proposition~\ref{wholomult2}.
    With this definition of multiplication, one can verify that
    \begin{equation*}
        R^f = 1 - \nabla_f U^f,
    \end{equation*}
    if $f$ is weakly holomorphic. This can be seen by using that this formula holds in the normalization,
    together with the fact that $U^{f'}$ has the standard extension property, SEP, i.e., that ${\bf 1}_{\{h = 0 \}} U^f = 0$
    for any tuple $h$ of strongly holomorphic functions not vanishing on any irreducible component of $Z$.
    This follows from that $U^{f'}$ is a principal value current, i.e., when $U^{f'}$ is written as a sum of
    push-forwards of elementary currents, the elementary currents contain no residue factors, and
    hence have the SEP.
\end{remark}

\begin{thm} \label{bmch}
    If $f = (f_1,\dots,f_p)$ is weakly holomorphic forming a complete
    intersection and $R^f = \mu \wedge e$, where $e = e_p\wedge \dots \wedge e_1$, then
\begin{equation*}
    \mu = \mu^f := \dbar \frac{1}{f_1}\wedge\cdots\wedge\dbar\frac{1}{f_p}.
\end{equation*}
\end{thm}

\begin{proof}
To begin with, we will assume that $f$ is strongly holomorphic. The proof will follow the same idea as the proof in the
smooth case in \cite{A2}, Theorem 3.1. Let
\begin{equation*}
    V = \frac{1}{f_p} e_p + \frac{1}{f_{p-1}}\dbar\frac{1}{f_p}\wedge e_p\wedge e_{p-1} + \cdots +
    \frac{1}{f_1}\dbar\frac{1}{f_2}\wedge\cdots\wedge\frac{1}{f_p}\wedge e_p\wedge\cdots\wedge e_1.
\end{equation*}
Then, by Proposition~\ref{commfimult}, $V$ satisfies
\begin{equation*}
    \nabla_f V = 1 - \dbar\frac{1}{f_1}\wedge\cdots\wedge\dbar\frac{1}{f_p} \wedge e.
\end{equation*}
Following the proof of Theorem 3.1 in \cite{A2}, locally, assume $Z \subseteq \Omega \subseteq \Cn$, $\omega$ is an arbitrary
neighborhood of $Z_f$ in $\Omega$ and $\chi$ is a smooth function with support on $\omega$ which is $\equiv 1$ in a neighborhood of $Z_f$.
Let $i : Z \to \Omega$ be the inclusion, and let $g = i^*\chi - i^*(\dbar \chi)\wedge u$. Then, since $\nabla_f u = 1$ on
$\supp \dbar\chi$, $\nabla_f g = 0$, and hence
\begin{align}  \label{nablafguv}
    \nabla_f( g\wedge(U^f-V) ) = g\wedge\nabla_f(U^f-V) = g_0(\mu^f-\mu)\wedge e = (\mu^f - \mu)\wedge e,
\end{align}
where $g_0 = \chi$ is the component of bidegree $(0,0)$ in $g$, which is $1$ in a neighborhood of $\supp (\mu^f-\mu)$.
A current $T$ is said to have the standard extension property, SEP, with respect to an analytic variety $W$ if for any
holomorphic function $h$ such that $h$ is not identically $0$ on any irreducible component of $W$, then $|h|^{2\lambda} T|_{\lambda = 0} = T$.
Since $\mu$ and $\mu^f$ are currents in $\PM(Z)$ of bidegree $(0,p)$, with support on $W = \{ f = 0 \}$,
$\mu$ and $\mu^f$ have the SEP, since if $h$ does not vanish on any irreducible component of $W$,
$\mu - |h|^{2\lambda}\mu|_{\lambda=0}$ has support on $W \cap \{ h = 0 \}$, which has codimension $\geq p+1$, and by
Proposition~\ref{pseudom0} it is $0$. Also, $\mu$ and $\mu^f$ are $\dbar$-closed and are annihilated by $\bar{I}_W$,
see Proposition~\ref{pseudom0}, so $i_*\mu, i_*\mu^f \in CH_W$,
where $CH_W$ denotes $\dbar$-closed $(0,\codim W)$-currents with support on $W$ satisfying the SEP.
By Lemma 3.3 in \cite{A2}, we know that a $\dbar$-closed current in $CH_W$ cannot be equal to
$\dbar \nu$, where $\nu$ can be chosen with support arbitrarily close to $W$, unless it is $0$. Hence, by looking at the
components of top degree in \eqref{nablafguv},
we have $i_*(\mu - \mu^f) = 0$, so $\mu = \mu^f$.

Now, if $f_i$ are weakly holomorphic, then the current $R^f$ will be the push-forward of the corresponding current $R^{\pi^* f}$,
where $\pi : Z' \to Z$ is the normalization of $Z$, and the same holds for the Coleff-Herrera product $\mu^f$. Hence,
equality holds in the normalization, and taking push-forward we get equality in the general case.
\end{proof}

\section{The transformation law}

With the Bochner-Martinelli type currents developed in the previous section, we will now prove the transformation law for Coleff-Herrera products
of weakly holomorphic functions. 

\begin{thm} \label{transformation}
    Assume that $f = (f_1,\cdots,f_p)$ and $g=(g_1,\cdots,g_p)$ are weakly holomorphic, defining complete intersections, and
    that there exists a matrix $A$ of weakly holomorphic
    functions such that $g = Af$. Then
    \begin{equation*}
        \dbar \frac{1}{f_1}\wedge\cdots\wedge \dbar \frac{1}{f_p} = (\det A)\dbar \frac{1}{g_1}\wedge\cdots\wedge\dbar\frac{1}{g_p}.
    \end{equation*}
\end{thm}

If $A$ is invertible, one can prove the transformation law with the help of Theorem~\ref{bmch} together
with the fact that the Bochner-Martinelli current is independent of the metric chosen to define $\sigma^f$
(here, in \eqref{sigmau}, $\sigma^f$ is defined with respect to the trivial metric on $E$), see \cite{A1}.
We will see that we can use a similar idea even in the case that $A$ is not invertible.
In \cite{Den} Denkowski proved the transformation law
for c-holomorphic functions based on a more direct approach.

To begin with, we assume that $f$, $g$ and $A$ are strongly holomorphic.
As in the previous section, we will identify $f$ and $g$ with sections of vector bundles,
however we will here identify them with sections of two different vector bundles.
Let $E$ and $E'$ be trivial holomorphic vector bundles over $Z$ with frames $e$ and $e'$, and make the identifications
$f = \sum f_i e_i^*$, $g = \sum g_i e_i'^*$ and $A \in \Hom(E',E)$ such that $g = f A$.

\begin{lma} \label{gdeltafprim}
    Let $\bigwedge A : \bigwedge E' \to \bigwedge E$ denote the linear extension of the mapping $(\bigwedge A) (v_1 \wedge \dots \wedge v_k) = Av_1\wedge \dots \wedge Av_k$.
    Then $\delta_f(\bigwedge A) = (\bigwedge A)\delta_g$.
\end{lma}

\begin{proof}
    Note first that $\delta_f A e_j' = g_j = \delta_{g} e_j'$. Hence, we have
    \begin{align*}
        &\delta_f(\bigwedge A)(e_{i_1}'\wedge\cdots\wedge e_{i_k}') = \delta_f( A e_{i_1}'\wedge\cdots\wedge A e_{i_k}') \\
        &= \sum (-1)^{j-1} A e_{i_1}'\wedge\cdots\wedge \delta_f( A e_{i_j}')\wedge\cdots\wedge Ae_{i_k}' \\
        &= \sum (-1)^{j-1} (\bigwedge A)( e_{i_1}'\wedge\cdots\wedge \delta_{g} e_{i_j}'\wedge\cdots\wedge e_{i_k}')
        = (\bigwedge A) \delta_{g} (e_{i_1}'\wedge\cdots\wedge e_{i_k}').
    \end{align*}
\end{proof}

To relate the currents $\mu^f$ and $\mu^g$, we will first derive a relation between the currents $U^f$ and $U^g$ as defined
by \eqref{UR}. 

\begin{lma} \label{uguf}
    If $f$ and  $g$ are strongly holomorphic and defining complete intersections, then there exists a current $R_1$ such that
    $U^f - (\bigwedge A) U^g = \nabla_f R_1$. 
\end{lma}

\begin{proof}
    Let $\sigma,u,\sigma'$ and $u'$ be the forms defined by \eqref{sigmau} corresponding to $f$ and $g$.
    Since $A$ is holomorphic, $(\bigwedge A) \dbar\sigma' = \dbar(A\sigma')$ outside of $\{ g = 0 \}$, and hence if
    we let $u_A' = \sum (A\sigma')\wedge (\dbar A\sigma')^{k-1}$,
    then $\nabla_f u_A' = 1$ outside of $\{ g = 0 \}$ by Lemma~\ref{gdeltafprim}. Thus, if $\Re \lambda \gg 0$,
    \begin{equation} \label{ugminusu}
        \nabla_f (|g|^{2\lambda} u_A'\wedge u) = |g|^{2\lambda} u - |g|^{2\lambda} u_A'
          - \dbar |g|^{2\lambda}\wedge u_A'\wedge u.
    \end{equation}
    We want to see that all the terms in \eqref{ugminusu} have current-valued analytic continuations to $\lambda = 0$.
    First, we note that since $\{ g = 0 \} \supseteq \{ f = 0 \}$, $|g|^{2\lambda} u |_{\lambda = 0} = U^f$ by Proposition~\ref{bmexist},
    and since $u_A' = (\bigwedge A) u'$ we get that $|g|^{2\lambda} u_A'|_{\lambda = 0} = (\bigwedge A) U^g$.
    Thus it remains to see that the left-hand side of \eqref{ugminusu} has an analytic continuation to $\lambda = 0$,
    and that the analytic continuation of the last term vanishes at $\lambda = 0$.
    To see that those terms have analytic continuations to $\lambda = 0$ is similar to showing the existence of the Bochner-Martinelli
    currents $U^f$ and $R^f$. If we recall briefly the proof of the existence of $U^f$ and $R^f$ in \cite{A1},
    the key step was that $\sigma\wedge (\dbar \sigma)^{k-1}$ is homogeneous with respect to $f$ in the sense
    that if $f = f_0 f'$, then $\sigma \wedge (\dbar \sigma)^{k-1} = (1/f_0^k) \sigma_0\wedge(\dbar \sigma_0)^{k-1}$,
    where $\sigma_0$ is smooth if $|f'| \neq 0$. 
    By blowing up along the ideals $(f_1,\dots,f_p)$ and $(g_1,\dots,g_p)$ followed by a resolution of singularities, see \cite{AHV},
    we can assume that locally $\pi^* f = f_0 h$ and $\pi^* g = g_0 g'$, where $h \neq 0$, $g' \neq 0$,
    and by a further resolution of singularities, we can assume that locally $f_0, g_0$ are monomials.
    Since $\{ g = 0 \} \supseteq \{ f = 0\}$, we get that $\{ g_0 = 0 \} \supseteq \{ f_0 = 0 \}$. 
    Thus, by the homogeneity of $\sigma' \wedge (\dbar\sigma')^{k-1}$ and $\sigma\wedge(\dbar \sigma)^{l-1}$ with respect to
    $f$ and $g$, we get, since $u_A' = (\bigwedge A) u'$, that $|g|^{2\lambda} u_A'\wedge u$ and $\dbar |g|^{2\lambda} \wedge u_A'\wedge u$
    acting on a test form $\varphi$ becomes finite sums of the form
    \begin{equation*}
        \int \frac{|v|^{2\lambda} |g_0|^{2\lambda}}{(g_0)^k f_0^l} \xi_{k,l}\wedge\pi^*\varphi \quad \text{and} \quad
        \int \frac{\dbar(|v|^{2\lambda} |g_0|^{2\lambda})}{(g_0)^k f_0^l}\wedge \xi_{k,l}\wedge\pi^*\varphi,
    \end{equation*}
    where $\xi_{k,l}$ are smooth $(0,k + l - 2)$-forms. Thus both have analytic continuations to $\lambda = 0$,
    and $R_2 := \dbar |g|^{2\lambda} \wedge u_A' \wedge u|_{\lambda = 0}$ has support on $\{ g = 0 \}$. Since $R_2 \in \PM(Z)$
    and consists of terms of bidegree $(0,k + l - 1)$, where $k + l \leq p$, with support on $\{ g = 0 \}$ which
    has codimension $p$, we get that $R_2 = 0$ by Proposition~\ref{pseudom0}. Thus, if we let $R_1 := |g|^{2\lambda} u_A'\wedge u|_{\lambda = 0}$,
    we get that $\nabla_f R_1 = U^f - (\bigwedge A) U^g$.
\end{proof}

Now we are ready to prove the transformation law.

\begin{proof} [Proof of Theorem \ref{transformation}]
    Assume first that $f,g$ and $A$ are strongly holomorphic, and make the same identifications as 
    after the statement of Theorem~\ref{transformation}.
    Since $(\bigwedge A)R^g = (\bigwedge A)(1-\nabla_g U^g) = 1-\nabla_f (\bigwedge A) U^g$ by Lemma~\ref{gdeltafprim},
    we get from Lemma~\ref{uguf} that
    \begin{equation*}
        (\bigwedge A) R^g - R^f = \nabla_f \left( (\bigwedge A) U^g - U^{f}\right) =
        \nabla_{f}^2 R_1 = 0,
    \end{equation*}
    so
    \begin{equation*}
        (\bigwedge A) R^g = R^f.
    \end{equation*}
    Thus, we get by Theorem~\ref{bmch} that
    \begin{equation*}
        (\bigwedge A) \left( \mu^g  \wedge e_p'\wedge\cdots\wedge e_1' \right) =
        \mu^f \wedge e_p\wedge\cdots\wedge e_1,
    \end{equation*}
    and since the left-hand side is equal to
    \begin{equation*}
        (\det A) \mu^g \wedge e_p\wedge\cdots\wedge e_1,
    \end{equation*}
    the transformation law follows. Now, if $f,g$ and $A$ are weakly holomorphic, the transformation law follows since
    equality must hold in the normalization because the pullback of $f$ and $g$ define complete intersections in the normalization.
    Hence, equality must hold also in $Z$ by taking push-forward.
\end{proof}

\section{The Poincar\'e-Lelong formula}

Let $f_1,\cdots,f_p$ be strongly holomorphic functions forming a complete intersection. The Poincar\'e-Lelong formula says
that
\begin{equation} \label{pl}
    \frac{1}{(2\pi i)^p}\dbar \frac{1}{f_1}\wedge\cdots\wedge\dbar \frac{1}{f_p}\wedge df_p\wedge\cdots\wedge df_1 = [Z_f] = \sum \alpha_i [V_i],
\end{equation}
where $V_i$ are the irreducible components of $Z_f$ and $[Z_f]$ is the integration current on $Z_f$ with multiplicities.
In case $p = \dim Z$ the multiplicity $\alpha_i$ at a point $x_i \in Z_f$ is given as the number of elements near $x_i$ of a generic fiber of $f$.
In case $p < \dim Z$ the multiplicity is given as the intersection multiplicity of $Z_f$ with $L$, where $L$ is a plane
of dimension $\dim Z - p$ transversal to $Z_f$. For a thorough discussion of the multiplicities see \cite{Ch},
and for a proof of the Poincar\'e-Lelong formula see Section 3.6 in \cite{CH}.

Now, if $f_i$ are weakly holomorphic functions defining a complete intersection,
we can give a relatively short proof that a formula similar to \eqref{pl} holds in $Z$.
In the strongly holomorphic case, assuming $Z\subseteq \Omega \subseteq \Cn$, $i_*[Z_f]$ can be seen either
as the intersection of the holomorphic chains $Z_{F_i}$ with $Z$, where $F_i$ are some holomorphic extensions of $f_i$ to $\Omega$,
or as a product of closed positive currents, see \cite{Ch}, that is
\begin{equation*}
    i_*[Z_f] = [Z_{F_1} \cdot \cdots\cdot Z_{F_p}\cdot Z] = [Z_{F_1}]\wedge\cdots\wedge[Z_{F_p}]\wedge[Z].
\end{equation*}
However, these types of products are in general only defined in case $Z_{F_1}\cap\cdots\cap Z_{F_p}\cap Z$ has codimension equal to
$\codim Z + \sum \codim Z_{F_i}$. Since zero sets of weakly holomorphic functions are
in general not zero sets of strongly holomorphic functions, as we saw in Example~\ref{zerosetnotholo}, we cannot expect to have a
similar interpretation for weakly holomorphic functions, since there are no natural counterparts to the holomorphic
$(n - 1)$-chains $Z_{F_i}$ or closed positive $(1,1)$-currents $[Z_{F_i}]$.

From now on, we assume that $f = (f_1,\dots,f_p)$ is weakly holomorphic defining a complete intersection.
Let $\pi : Z' \to Z$ be the normalization of $Z$, so that in particular, $\pi$ is a finite proper holomorphic map.
Since $f' = \pi^* f$ forms a complete intersection,
\eqref{pl} holds for $f'$ in the normalization. Note that, if $V_i$ are the irreducible components of $Z_{f'}$, then
$W_i := \pi(V_i)$ are irreducible in $Z$.
If $f : V \to W$ is a branched holomorphic cover with exceptional set $E$, we say that $f$ is a \emph{*-covering} if $W\setminus E$
is a connected manifold. In particular, this means that the sheet-number of $f$ is constant outside the exceptional set.
By the Andreotti-Stoll theorem, see \cite{Loj}, if $f : V \to W$ is a finite proper holomorphic map,
$V$ has constant dimension and $W$ is irreducible, then $f$ is a *-covering.
If $V \subset Z'$ is an irreducible component of $Z_{f'}$ and we consider $\pi|_V : V \to W$, where $W = \pi(V)$, it is a finite
proper holomorphic map satisfying the conditions required for the Andreotti-Stoll theorem. Hence, there exists an
integer $k$ such that $\pi|_V$ is a $k$-sheeted finite branched holomorphic covering.
Thus $\pi_* \alpha [V] = k\alpha [W]$.
For $f = (f_1,\cdots,f_p)$ a weakly holomorphic mapping forming a complete intersection, we define the left-hand side of \eqref{pl}
as the push-forward of the corresponding current in the normalization. Thus, since we have by \eqref{pl} that
\begin{equation*}
    \frac{1}{(2\pi i)^p}\dbar \frac{1}{f_1}\wedge\cdots\wedge\dbar \frac{1}{f_p}\wedge df_p\wedge\cdots\wedge df_1 =
    \pi_*[Z_{f'}],
\end{equation*}
we have proved the following.
\begin{thm} \label{plthm}
    Let $f = (f_1,\cdots,f_p)$ be a weakly holomorphic mapping forming a complete intersection. Then
\begin{equation}
    \frac{1}{(2\pi i)^p}\dbar \frac{1}{f_1}\wedge\cdots\wedge\dbar \frac{1}{f_p}\wedge df_p\wedge\cdots\wedge df_1 = \sum \beta_i [W_i]
\end{equation}
where $\beta_i \in \mathbb{N}$ and $W_i$ are the irreducible components of $W = Z_f$.
More explicitly, if $[Z_{f'}] = \sum \alpha_i [V_i]$ and say $V_{i_1},\cdots,V_{i_k}$ are the sets $V_j$ such that
$\pi(V_j) = W_i$, then $\beta_i = \sum k_{i_j}\alpha_{i_j}$, where $k_j$ is the number of elements in a generic fiber of
$\pi|_{V_j}$.
\end{thm}

\begin{remark}
    In \cite{Den} Denkowski proves the Poincar\'e-Lelong formula for $f = (f_1,\dots,f_p) \in \Oc^{\oplus p}(Z)$
    (based on his construction on $\Gamma_f$, however as for the Coleff-Herrera product in Proposition~\ref{denequal}
    our definition coincides with his).
    In that case, it gives a different interpretation of the multiplicities as the intersection cycle
    \begin{equation*}
        \frac{1}{(2\pi i)^p}\dbar\frac{1}{f_1}\wedge \dots \wedge \dbar\frac{1}{f_p}\wedge df_p\wedge\dots\wedge df_1
        = \pi_* ([\Gamma_f]\cdot[Z\times \{ 0 \}]),
    \end{equation*}
    where $\pi : Z\times \C^p \to Z$ is the projection.
\end{remark}

Note that if $f$ is weakly holomorphic, since $f$ is in general not smooth on $Z_\sing$, $df$ is not in general defined on all $Z$
(although its pullback to the normalization has a smooth extension to all of $Z'$) so, as for multiplication with weakly holomorphic functions
in Example~\ref{exmultcholo}, it might for example happen that $\dbar (1/f) = 0$ while $\dbar (1/f) \wedge df \neq 0$.
For example, if $Z = \{ z^3 = w^2 \}$, $\pi(t) = (t^2,t^3)$ and $f = w/z \in \Ow(Z)$, that is $\pi^* f = t$, then $\dbar (1/f) = 0$ while
$\dbar (1/f) \wedge df = 2\pi i [0]$, as expected, since $Z_f = \{ 0 \}$.

\section*{Acknowledgments}

I would like to thank Mats Andersson for suggesting the original idea of this article,
and for interesting discussions about the topic.

\end{document}